\documentclass[11pt, reqno, oneside]{amsart}

\usepackage{geometry}
\usepackage{float}
\usepackage{mathrsfs}
\usepackage[active]{srcltx}
\usepackage{mathrsfs}
\usepackage{array}
\usepackage{multirow}
\usepackage{float}
\usepackage{caption}
\usepackage{subfig}
\usepackage{yfonts}
\usepackage[final]{pdfpages}
\DeclareMathAlphabet{\mathpzc}{OT1}{pzc}{m}{rm}

\usepackage{color}
\usepackage[unicode]{hyperref}
\hypersetup{
	colorlinks = true,%
	citecolor = [rgb]{0.0,0.0,0.9},
	filecolor=black,%
	linkcolor = [rgb]{0.65,0.0,0.0},%
	anchorcolor = red,
	pagecolor = red,
	urlcolor= [rgb]{0.65,0.0,0.0}
}

\numberwithin{equation}{section}

\usepackage{graphicx}

\newtheorem{example}{ Example}[section]
\newtheorem{proposition}{Proposition}[section]
\newtheorem{theorem}{Theorem}[section]

\newtheorem{definition}{Definition}[section]
\newtheorem{remark}{Remark}[section]

\newcommand{\ds}{\displaystyle}

\title[TITLE ~~~\today]{Differential inclusions involving oscillatory terms}
\author[Authors]{Alexandru Krist\'aly, Ildik\'o I. Mezei, K\'aroly Szil\'ak}
\thanks{A. Krist\'aly was supported by the National Research, Development and Innovation Office, Hungary, K-18, project no. 127926.}
\thanks{Dedicated to Professor Gheorghe Moro\c sanu on the occasion of his 70th birthday}
\address{Institute of Applied Mathematics\\
	\'Obuda University, 1034 Budapest, Hungary \& Department of Economics\\
	Babe\c s-Bolyai University\\
	400591 Cluj-Napoca, Romania 
	}
\email{kristaly.alexandru@nik.uni-obuda.hu; alex.kristaly@econ.ubbcluj.ro}

\address{Department of Mathematics and Informatics
	Babe\c s-Bolyai University
	400084 Cluj-Napoca, Romania}
\email{ildiko.mezei@math.ubbcluj.ro}

\address{Institute of Applied Mathematics\\
	\'Obuda University\\
	1034 Budapest, Hungary}
\email{karoly.szilak@gmail.com}


\begin{document}

	\begin{abstract}
	 Motivated by mechanical problems where external forces are non-smooth,  we consider the differential inclusion problem  
	 \[ \   \left\{ \begin{array}{lll}
	 -\Delta u(x)\in \partial F(u(x))+\lambda \partial G(u(x))& {\rm in} &   \Omega; \\
	 u\geq 0, &\mbox{in} &   \Omega;\\
	 u= 0, &{\rm on}& \partial \Omega,
	 \end{array}\right. \eqno{({\mathcal D}_\lambda)}\]
	 where $\Omega \subset {\mathbb R}^n$ is a bounded open domain, and $\partial F$ and $\partial G$ stand for the  generalized gradients of the locally Lipschitz functions $F$ and $G$.   In this paper we provide a quite complete picture on the number of solutions of $({\mathcal D}_\lambda)$ whenever $\partial F$ oscillates near the origin/infinity and $\partial G$ is a generic perturbation of order $p>0$ at the origin/infinity, respectively. Our results extend in several aspects those of Krist\'aly and Moro\c sanu [\textit{J. Math. Pures Appl}., 2010]. 
	\end{abstract}
	\maketitle

	\vspace{-1.1cm}
	\section{Introduction}

	
	We consider the model Dirichlet problem 
	\[ \   \left\{ \begin{array}{lll}
	-\Delta u(x)= f(u(x))& {\rm in} &   \Omega; \\
	u\geq 0, &\mbox{in} &   \Omega;\\
	u= 0, &{\rm on}& \partial \Omega,
	\end{array}\right. \eqno{({P_0})}\]
	where $\Delta$ is the usual Laplace operator, $\Omega \subset {\mathbb R}^n$ is a bounded open domain $(n\geq 2)$, and $f:\mathbb R\to \mathbb R$ is a continuous function  verifying certain growth conditions at the origin and infinity. Usually, such a problem is studied on the Sobolev space $H_0^1(\Omega)$ and weak solutions of $(P_0)$ become classical/strong solutions whenever $f$ has further regularity. There are several approaches to treat problem $(P_0)$, mainly depending on the behavior of the function $f$. When $f$ is superlinear and subcritical at infinity (and superlinear at the origin), the seminal paper of Ambrosetti and Rabinowitz \cite{AR} guarantees the existence of at least a nontrivial solution of $(P_0)$ by using variational methods.  
	An important extension of $(P_0)$ is its \textit{perturbation}, i.e., 
	\[ \   \left\{ \begin{array}{lll}
	-\Delta u(x)= f(u(x))+\lambda g(u(x)) & {\rm in} &   \Omega; \\
	u\geq 0, &\mbox{in} &   \Omega;\\
	u= 0, &{\rm on}& \partial \Omega,
	\end{array}\right. \eqno{({P_\lambda})}\]
	where $g:\mathbb R\to \mathbb R$ is another continuous function which is going to compete with the original function $f$. When both functions $f$ and $g$ are of \textit{polynomial type} of sub- and super-unit degree, --  the right hand side being called as a concave-convex nonlinearity -- the existence of at least one or two nontrivial solutions of $(P_\lambda)$ is guaranteed, depending on the range of $\lambda>0$, see e.g. Ambrosetti, Brezis and Cerami \cite{ABC}, Autuori and Pucci \cite{AP}, de Figueiredo,  Gossez and Ubilla \cite{Fig}. In these papers variational arguments, sub- and super-solution methods as well as fixed point arguments are employed. 
	
	Another important class of problems of the type $(P_\lambda)$ is studied whenever $f$ has a certain \textit{oscillation} (near the origin or at infinity) and $g$ is a \textit{perturbation}. Although oscillatory functions seemingly call forth the existence of infinitely many  solutions, it turns out that 'too classical' oscillatory functions do not have such a feature. Indeed, when $f(s)=c\sin s$ and $g=0$, with $c>0$ small enough,  a simple use of the Poincar\'e inequality implies that problem $(P_\lambda)$ has only the zero solution. However, when $f$ \textit{strongly} oscillates, problem $(P_0)$ has indeed infinitely many different solutions; see e.g. Omari and Zanolin \cite{OZ}, Saint Raymond \cite{Ray}. Furthermore, if $g(s)=s^p$ $(s>0)$, a novel competition phenomena has been described for  $({P_\lambda})$ by Krist\'aly and Moro\c sanu \cite{KM}. 	We notice that several extensions of \cite{KM} can be found in the literature, see e.g. Ambrosio, D'Onofrio and Molica Bisci \cite{ADMB} and Molica Bisci and Pizzimenti \cite{MBP} for nonlocal fractional Laplacians; Molica Bisci, R\u adulescu and Servadei \cite{MBRS} for general operators in divergence form; M\u alin and R\u adulescu \cite{MR} for difference equations. We emphasize that in the aforementioned papers the perturbations are either zero or have a (smooth) polynomial form.

	In mechanical applications, however, the perturbation may occur in a \textit{discontinuous} manner as a non-regular external force, see e.g. the gluing force in von K\'arm\'an laminated plates, cf. Bocea,   Panagiotopoulos and   R\u adulescu \cite{Bocea},   Motreanu and  Panagiotopoulos \cite{minimax} and Panagiotopoulos \cite{Pan}. In order to give a reasonable reformulation of problem $({P_\lambda})$ in such a non-regular setting, 
 the idea is to  'fill the gaps' of the discontinuities, considering instead of the discontinuous nonlinearity a \textit{set-valued map} appearing as the generalized gradient of a locally Lipschitz function. In this way, we deal with an \textit{elliptic differential inclusion} problem  rather than an elliptic differential equation, see e.g.  Chang \cite{Chang}, Gazzolla and R\u adulescu \cite{GR} and Krist\'aly \cite{Kris-JDE}; this problem can be formulated generically as 
	\[ \   \left\{ \begin{array}{lll}
	-\Delta u(x)\in \partial F(u(x))+\lambda \partial G(u(x))& {\rm in} &   \Omega; \\
	u\geq 0, &\mbox{in} &   \Omega;\\
	u= 0, &{\rm on}& \partial \Omega,
	\end{array}\right. \eqno{({\mathcal D}_\lambda)}\]
	where  $F$ and $G$ are both nonsmooth, locally Lipschitz functions having various growths, while $\partial F$ and $\partial G$ stand for the  generalized gradients of $F$ and $G$, respectively.
	
	The main purpose of the present paper is to extend the main results  of Krist\'aly and Moro\c sanu \cite{KM} in two directions: 
\begin{itemize}
	\item[(a)] to allow the presence of nonsmooth nonlinear terms -- reformulated into the inclusion $({\mathcal D}_\lambda)$ -- which are more suitable from mechanical point of view (mostly due to the perturbation term $G$, although we allow non-smoothness for the oscillatory term $F$ as well); 
	
	\item[(b)] to consider a generic $p$-order perturbation $\partial G$  at the origin/infinity, not necessarily of polynomial growth as in \cite{KM}, $p>0.$
\end{itemize}

	In the present paper we study the inclusion $({\mathcal D}_\lambda)$ in two different settings, i.e., we analyze the number of distinct solutions of $({\mathcal D}_\lambda)$ whenever  $\partial F$ oscillates near the origin/infinity and  $\partial G$ is  of order $p>0$ near the origin/infinity. 	Roughly speaking, when $\partial F$ \textit{oscillates near the origin} and \textit{$\partial G$ is of order $p>0$ at the origin},  we prove  that the number of distinct, nontrivial  solutions  of $({\mathcal D}_\lambda)$  is
	\begin{itemize}
		\item infinitely many whenever  $p>1$ ($\lambda\geq 0 $ is arbitrary) or  $p=1$ and $\lambda$ is small enough (see Theorem \ref{elso-tetel}); 
		\item at least (a prescribed number) $k\in \mathbb N$  whenever $0<p<1$ and $\lambda$ is small enough (see Theorem \ref{masodik-tetel}).
	\end{itemize}
	As we can observe, in the first case, the term $\partial G(s)\sim s^p$ as $s\to 0^+$ with $p>1$ has no effect on the number of solutions (i.e., the oscillatory  term is the leading one), while in the second case, the situation changes dramatically, i.e., $\partial G$ has a 'truth' competition with respect to the  oscillatory term $\partial F$. 
	
	We can state a very similar result as above whenever \textit{$\partial F$ oscillates at infinity} and  \textit{$\partial G$ is of order $p>0$ at infinity} by proving that the number of distinct, nontrivial solutions of the  differential inclusion $({\mathcal D}_\lambda)$ is
	\begin{itemize}
		\item infinitely many whenever  $p<1$ ($\lambda\geq 0$ is arbitrary) or  $p=1$ and $\lambda$ is small enough (see Theorem \ref{harmadik tetel}); 
		\item at least (a prescribed number) $k\in \mathbb N$  whenever $p>1$ and $\lambda$ is small enough (see Theorem \ref{negyedik tetel}).
	\end{itemize}
	Contrary to the competition at the origin, in the first case the term $\partial G(s)\sim s^p$ as $s\to \infty$ with $p<1$ has no effect on the number of solutions (i.e., the oscillatory  term is the leading one), while in the second case, the perturbation term $\partial G$ competes with the oscillator function $\partial F$.

	We admit that the line of the proofs is conceptually similar to that of Krist\'aly and Moro\c sanu \cite{KM}; however, the presence of the nonsmooth terms $\partial F$ and $\partial G$ requires a deep argumentation by fully exploring the nonsmooth calculus of locally Lipschitz functions in the sense of Clarke \cite{Clarke}. In addition, the presence of the generic $p$-order perturbation $\partial G$ needs a special attention with respect to \cite{KM}; in particular, the  $p$-order growth of $\partial G$ is new even in smooth settings. 

The organization of the present paper is the following.  In Section \ref{Main theorems} we state our main assumptions and results, providing also some examples of functions fulfilling the assumptions. Section \ref{generic result} contains a generic localization theorem for differential inclusions, while Sections \ref{proof of theorems with osc near orig} and \ref{proof of theorems with osc at infinity} are devoted to the proof of our main results. In Section \ref{conc-section} we formulate some concluding remarks, while in the Appendix (Section \ref{appendix}) we  collect those notions and results on locally Lipschitz functions that are used throughout our arguments.

\section{Main theorems}\label{Main theorems}

	Let $F,G:\mathbb R_+\to \mathbb R$ be locally Lipschitz functions and as usual, let us denote by $\partial F$ and $\partial G$ their generalized gradients in the sense of Clarke (see the Appendix). Hereafter, $\mathbb R_+=[0,\infty).$ Let $p>0$, $\lambda\geq 0$ and $\Omega\subset \mathbb R^n$ be a bounded open domain, and  consider the elliptic differential inclusion problem 
	 \[ \   \left\{ \begin{array}{lll}
	 -\Delta u(x)\in \partial F(u(x))+\lambda \partial G(u(x))& {\rm in} &   \Omega; \\
	  u\geq 0 &\mbox{in} &   \Omega;\\
	 u= 0, &{\rm on}& \partial \Omega.
	 \end{array}\right. \eqno{({\mathcal D}_\lambda)}\]

We distinguish the cases when $\partial F$ oscillates near the \textit{origin} or at \textit{infinity}.

\subsection{Oscillation near the origin}
	We assume:
	\begin{enumerate}
		\item[$(F^0_0)$]  $F(0)=0$;
		\item[$(F^0_1)$]
			$-\infty<\liminf_{s\to
			0^+}\frac{F(s)}{s^2}; \  \limsup_{s\to
			0^+}\frac{F(s)}{s^2}=+\infty$; 
		\item[$(F^0_2)$]
			$l_0:=\liminf_{s\to 0^+}\frac{\max\{\xi:\xi\in \partial F(s)\}}{s}<0.$ 
		\item[$(G^0_0)$] $G(0)=0$;
		\item[$(G^0_1)$]
			There exist $p>0$ and $\underline c,\overline c\in \mathbb R$ such that 
			$$\underline c=\liminf_{s\to 0^+}\frac{\min\{\xi:\xi\in \partial G(s)\}}{s^p}\leq \limsup_{s\to 0^+}\frac{\max\{\xi:\xi\in \partial G(s)\}}{s^p}=\overline c.$$
	\end{enumerate}
\begin{remark}\label{rem-1} \rm  Hypotheses $(F^0_1)$ and $(F^0_2)$
	imply a strong oscillatory behavior of $\partial F$ near the origin. Moreover, it turns out that $0\in \partial F(0)$; indeed, if we assume the contrary, by the upper semicontinuity of $\partial F$ we also have that $0\notin \partial F(s)$ for every small $s>0$. Thus, by $(F^0_2)$ we have that $\partial F(s)\subset (-\infty,0]$ for these values of $s>0$. By using $(F^0_0)$ and Lebourg's mean value theorem (see Proposition \ref{lebourg kozepertek tetel} in the Appendix), it follows that $F(s)=F(s)-F(0)=\xi s\leq 0$ for some $\xi\in \partial F(\theta s)\subset (-\infty,0]$ with $\theta\in (0,1)$. The latter inequality contradicts the second assumption from $(F^0_1)$. 
	Similarly, one obtains that $0 \in \partial G(0)$ by exploring $(G_0^0)$ and $(G_1^0)$, respectively.

	
	In conclusion, since $0\in \partial F(0)$ and $0\in \partial G(0)$, it turns out that $0\in H_0^1(\Omega)$ is a solution of the differential inclusion $({\mathcal D}_\lambda)$. Clearly, we are interested in nonzero solutions of $({\mathcal D}_\lambda)$. 
\end{remark}

\begin{figure}[H]
	\centering
	{{\includegraphics[width=7cm]{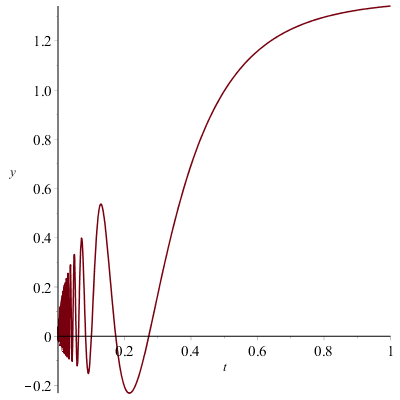} }}%
	{{\includegraphics[width=7cm]{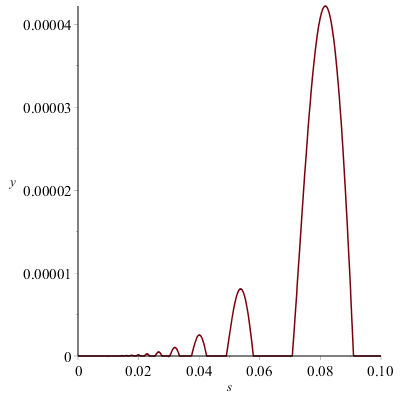} }}%
	\caption{Graphs of $f_0$ and $G_0$ around the origin, respectively.}%
	\label{fig. 2-odrendu origoban}%
\end{figure}

\begin{example}\rm 
	Let us consider $F_0(s)=\ds\int_0^s f_0(t)$, $s\geq 0$, where $f_0(t)=\sqrt{t}(\frac{1}{2}+\sin t^{-1})$, $t>0$ and $f_0(0)=0$, or some of its jumping variants. One can prove that $\partial F_0=f_0$ verifies the assumptions $(F^0_0)-(F^0_2)$.
	For a fixed  $p>0,$ let $G_0(s)=\ln (1+s^{p+2})\max\{0,\cos s^{-1}\}$, $s>0$ and $G_0(0)=0.$  It is clear that $G_0$ is not of class $C^1$ and verifies $(G_1^0)$ with $\underline c=-1$ and $\overline c=1$, respectively; see Figure \ref{fig. 2-odrendu origoban} representing both $f_0$ and $G_0$ (for $p=2$).  
	
\end{example}

	In the sequel, we provide a quite complete picture about the competition concerning the terms $s\mapsto \partial F(s)$ and $s\mapsto \partial G(s)$, respectively. First, we are going to show that when $p\geq 1$ then the 'leading' term is the oscillatory function $\partial F$; roughly speaking, one can say that the effect of $s\mapsto \partial G(s)$ is  negligible in this competition. More precisely, we prove the following result.

\begin{theorem}\label{elso-tetel} {\rm (Case $ p \geq 1$)} Assume that $ p \geq 1$ and the locally Lipschitz functions $F,G:\mathbb R_+\to \mathbb R$ satisfy $(F^0_0)-(F^0_2)$ and $(G^0_0)-(G^0_1).$ If 
	
	{\rm  (i)} either $ p =1$ and $\lambda \overline c<-l_0$ $($with $\lambda\geq 0$$)$,
	
	{\rm (ii)} or $ p >1$ and $\lambda\geq 0$ is
	arbitrary,\\
	then the differential inclusion problem $({\mathcal D}_\lambda)$ admits a sequence $\{u_i\}_{i}\subset H_0^1(\Omega)$ of distinct weak
	 solutions such that
	\begin{equation}\label{nulla-osc-hatarertek}
	\lim_{i\to
		\infty}\|u_i\|_{H_0^1}=\lim_{i\to \infty}\|u_i\|_{L^\infty}=0.
	\end{equation}
\end{theorem}

In the case when $p<1$, the perturbation term $\partial G$ may compete with the oscillatory function $\partial F$; namely, we have:

\begin{theorem}\label{masodik-tetel} {\rm (Case $0< p < 1$)} Assume $0< p< 1$ and that the locally Lipschitz functions $F,G:\mathbb R_+\to \mathbb R$ satisfy $(F^0_0)-(F^0_2)$ and $(G^0_0)-(G^0_1).$ Then, for every $k\in \mathbb N$, there
	exists $\lambda_k>0$ such that the differential inclusion $(\mathcal D_\lambda)$ has at least
	$k$ distinct weak solutions
	$\{u_{1,\lambda},...,u_{k,\lambda}\}\subset H_0^1(\Omega)$
	whenever $\lambda\in [0,\lambda_k].$ Moreover,
	\begin{equation}\label{becsles-i-vel}
	\|u_{i,\lambda}\|_{H_0^1}<i^{-1}\  \ and\ \
	\|u_{i,\lambda}\|_{L^\infty}<i^{-1}\ \  {for\ any}\
	i=\overline{1,k};\ \lambda\in [0,\lambda_k].
	\end{equation}
\end{theorem}

\vspace{0.5cm}

\subsection{Oscillation at infinity}
Let assume:
\begin{enumerate}
	\item[$(F^{\infty}_0)$]  $F(0)=0$;
	\item[$(F^{\infty}_1)$]
	$-\infty<\liminf_{s\to
		{\infty}}\frac{F(s)}{s^2}; \  \limsup_{s\to
		{\infty}}\frac{F(s)}{s^2}=+\infty$; 
	\item[$(F^{\infty}_2)$]
	$l_{\infty}:=\liminf_{s\to {\infty}}\frac{\max\{\xi:\xi\in \partial F(s)\}}{s}<0.$ 
	\item[$(G^{\infty}_0)$] $G(0)=0;$
	\item[$(G^{\infty}_1)$]
	There exist $p>0$ and $\underline c,\overline c\in \mathbb R$ such that 
	$$\underline c=\liminf_{s\to {\infty}}\frac{\min\{\xi:\xi\in \partial G(s)\}}{s^p}\leq \limsup_{s\to {\infty}}\frac{\max\{\xi:\xi\in \partial G(s)\}}{s^p}=\overline c.$$
\end{enumerate}
\begin{remark}  \rm
	Hypotheses $(F^\infty_1)$ and $(F^\infty_2)$
imply a strong oscillatory behavior of the set-valued map $\partial F$ at infinity.
	
\end{remark}
\begin{example}\rm 
	We consider $F_{\infty}(s)=\ds\int_0^s f_{\infty}(t)$, $s\geq 0$, where $f_{\infty}(t)=\sqrt{t}(\frac{1}{2}+\sin t)$, $t\geq 0$, or some of its jumping variants;  one has that $F_{\infty}$ verifies the assumptions $(F^{\infty}_0)-(F^{\infty}_2)$. 
	For a fixed $p>0$, let $G_{\infty}(s)=s^p\max\{0,\sin s\}$, $s\geq 0$;  it is clear that $G_{\infty}$ is a typically locally Lipschitz function on $[0,\infty)$ (not being of class $C^1$) and verifies $(G_1^{\infty})$ with $\underline c=-1$ and $\overline c=1$; see Figure \ref{fig. 2-odrendu inf}  representing both $f_\infty$ and $G_\infty$ (for $p=2$), respectively.
	
	\begin{figure}[H]
		\centering
		{{\includegraphics[width=7cm]{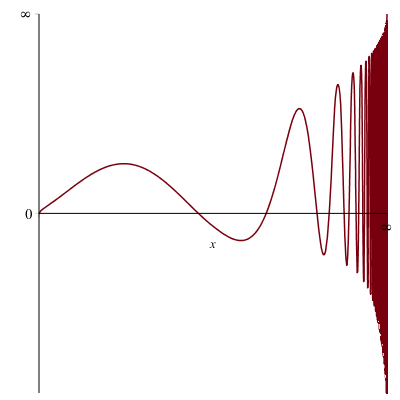} }}%
		{{\includegraphics[width=7cm]{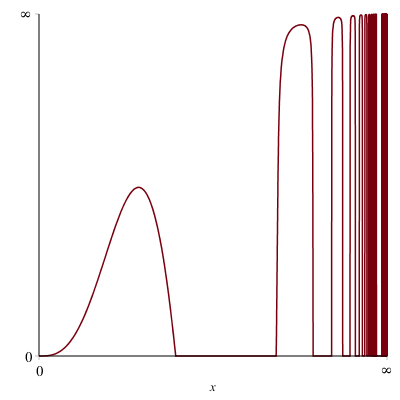} }}%
		\caption{Graphs of $f_\infty$ and $G_\infty$ at infinity, respectively.}%
		\label{fig. 2-odrendu inf}%
	\end{figure}

%
%
\end{example}

	In the sequel, we investigate the competition at infinity concerning the terms $s\mapsto \partial F(s)$ and $s\mapsto \partial G(s)$, respectively. 	First, we  show that when $p\leq 1$ then the 'leading' term is the oscillatory function $F$, i.e.,  the effect of $s\mapsto \partial G(s)$ is  negligible. More precisely, we prove the following result:

\begin{theorem} {\rm (Case $p \leq 1$)} \label{harmadik tetel} 
	Assume that $p\leq 1$ and the locally Lipschitz functions $F,G:\mathbb R_+\to \mathbb R$ satisfy $(F^{\infty}_0)-(F^{\infty}_2)$ and $(G^{\infty}_0)-(G^{\infty}_1)$. 
	 If 

	{\rm  (i)} either $ p =1$ and $\lambda \overline c \leq -l_0$ $($with $\lambda\geq 0$$),$

	{\rm (ii)} or $ p <1$ and $\lambda \geq 0$ is arbitrary,\\
	then the differential inclusion $({\mathcal D}_\lambda)$ admits a sequence $\{u_i\}_{i}\subset H_0^1(\Omega)$ of distinct weak
 solutions such that
	\begin{equation}\label{nulla-osc-hatarertek-infty}
		\lim_{i\to
		\infty}\|u_i^{\infty}\|_{L^{\infty}}=\infty.
	\end{equation}
\end{theorem}

\begin{remark} \rm
Let $2^*$ be the usual critical Sobolev exponent. In addition to (\ref{nulla-osc-hatarertek-infty}), we also have $\lim_{i\to \infty}\|u_i^{\infty}\|_{H_0^1}=\infty$ whenever
	\begin{equation}\label{inftycondF}
	\sup_{s\in [0, \infty)}\frac{\max\{|\xi|: \xi \in \partial F(s)\}}{1+s^{2^*-1}} < \infty.
	\end{equation}
\end{remark}

	\rm In the case when $p>1$, it turns out that the perturbation term $\partial G$ may compete with the oscillatory function $\partial F$; more precisely, we have:

\begin{theorem}{\rm (Case $p > 1$)} \label{negyedik tetel} Assume that $p > 1$ and the locally Lipschitz functions $F,G:\mathbb R_+\to \mathbb R$ satisfy $(F^{\infty}_0)-(F^{\infty}_2)$ and $(G^{\infty}_0)-(G^{\infty}_1)$.
Then, for every $k\in \mathbb N$, there
exists $\lambda_k^{\infty}>0$ such that the differential inclusion $(\mathcal D_\lambda)$ has at least
$k$ distinct weak solutions
$\{u_{1,\lambda},...,u_{k,\lambda}\}\subset H_0^1(\Omega)$
whenever $\lambda\in [0,\lambda_k^{\infty}].$ Moreover,
\begin{equation}\label{becsles-i-vel-infty}
\|u_{i,\lambda}\|_{L^\infty}>
i-1\ \  {for\ any}\
i=\overline{1,k};\ \lambda\in [0,\lambda_k^{\infty}].
\end{equation}
\end{theorem}

\begin{remark}\label{estimation_Tétel4}
	\rm If  (\ref{inftycondF}) holds and $p \leq 2^{*}-1$ in Theorem \ref{inftycondF}, then we have in addition that 
	$$\|{u_{i,\lambda}^{\infty}}\|_{H_0^1}>i-1\ \  {\rm for\ any}\
	i=\overline{1,k};\ \lambda\in [0,\lambda_k^{\infty}].$$
\end{remark}

\section{Localization: a generic result}\label{generic result} We consider the following differential inclusion 
problem { $$\left\{
	\begin{array}{lll}
	-\triangle u(x) + ku(x) \in   \partial A(u(x)),\ \ u(x)\geq 0 \ \ &  &
	x\in \Omega,
	\\ u(x)= 0 & &
	x\in \partial \Omega,
	\end{array}\right.
	\eqno{({\rm D}_A^k)}$$}
where $k>0$ and 
\begin{enumerate}
	\item[(H$_A^1$):] $A: [0,\infty)\to \mathbb R$ is a locally Lipschitz function with  $A(0)=0$, and there is $M_A>0$
	such that $$\max \{|\partial A(s)|\}:=\max\{|\xi|:\xi\in \partial A(s)\}\leq M_A$$ for every $s\geq 0$;
	\item[(H$_A^2$):]  there are $0<\delta<\eta$ such that $\max\{\xi:\xi\in \partial A(s)\}\leq 0$ for every 
 $s\in [\delta,\eta]$.
\end{enumerate}

For simplicity, we extend the function $A$ by $A(s)=0$ for $s\leq 0;$ the extended function is locally Lipschitz on the whole $\mathbb R.$ The natural energy functional $\mathcal
T:H_0^1(\Omega) \to \mathbb R$ associated with the differential inclusion problem $({\rm
	D}_A^k)$ is defined by $$\mathcal
T(u)=\frac{1}{2}\|u\|^2_{H^1_0}+\frac{k}{2}\int_{\Omega}u^2dx-\int_{\Omega}A(u(x))dx.\
$$

	The energy functional $\mathcal T$ is well defined and locally Lipschitz on $H_0^1(\Omega)$, while its critical points in the sense of Chang (see Definition \ref{kritpont-chang} in the Appendix) are precisely the weak solutions of the differential inclusion problem 
	$$\left\{
	\begin{array}{lll}
	-\triangle u(x) + ku(x) \in   \partial A(u(x)), \ \ &  &
	x\in \Omega,
	\\ u(x)= 0 & &
	x\in \partial \Omega;
	\end{array}\right.\eqno{({\rm D}_A^{k,0})}
	$$
	note that at this stage we have no information on the sign of $u$.
	
	Indeed, if $0\in \partial \mathcal T(u)$, then 
for every $v\in H_0^1(\Omega)$ we have
$$\int_{\Omega}\nabla u(x)\nabla v(x)dx-k\int_{\Omega}u(x) v(x){dx}-\int_{\Omega}\xi_x(x)v(x)dx=0,$$
where $\xi_x\in \partial A(u(x))$ a.e. $x\in \Omega$, see e.g. Motreanu and Panagiotopoulos \cite{minimax}. By using the divergence theorem for the first term at the left hand side (and exploring the Dirichlet boundary condition), we  obtain that
$$\int_{\Omega}\nabla u(x)\nabla v(x)dx=-\int_{\Omega}{\rm div}(\nabla u(x)) v(x)dx=-\int_{\Omega}\Delta u(x) v(x)dx.$$
Accordingly, we have that 
$$-\int_{\Omega}\Delta u(x) v(x)dx + k\int_{\Omega}u(x)v(x)=\int_{\Omega}\xi_x v(x){dx}$$
for every test function $v\in H_0^1(\Omega)$ which means that $-\Delta u(x) + ku(x) \in \partial A(u(x))$ in the weak sense in $\Omega$, as claimed before.\\

Let us consider the number $\eta\in \mathbb R$ from
(H$_A^2$) and the set $$W^\eta= \{u\in
H_0^1(\Omega):\|u\|_{L^\infty}\leq \eta\}. $$

Our localization result reads as follows (see  \cite[Theorem 2.1]{KM} for its smooth form): 

\begin{theorem}\label{segedtetel eredmeny oszcillacio nullaban} Let $k>0$ and assume that hypotheses  {\rm (H$_A^1$)} and {\rm (H$_A^2$)}
	hold. Then
	\begin{enumerate}
		\item[\rm (i)] the energy functional $\mathcal T$ is bounded from below
		on $W^\eta$ and its infimum is attained at some $\tilde u\in
		W^\eta$$;$
		
		\item[\rm (ii)] $\tilde u(x)\in [0,\delta]$ for a.e. $x\in
		\Omega$$;$
		
		\item[\rm (iii)]  $\tilde u$ is a weak solution of the differential inclusion $({\rm D}_A^{k}).$
	\end{enumerate}
\end{theorem}

{\it Proof.} The  proof is similar to that of Krist\'aly and Moro\c sanu \cite{KM}; for completeness, we provide its main steps. 

 $\rm (i)$ Due to (H$_A^1$), it is clear that the energy functional $\mathcal T$ is
 bounded from below on $H_0^1(\Omega)$.  Moreover, due to the compactness of the embedding $H_0^1(\Omega)\subset L^q(\Omega)$, $q\in [2,2^*)$, it turns out that $\mathcal T$ is sequentially weak lower semi-continuous on $H_0^1(\Omega)$. In addition, the set $W^\eta$ is weakly closed, being convex and
closed in $H_0^1(\Omega)$. Thus, there is  $\tilde u\in W^\eta$ which is a minimum
point of $\mathcal T$ on the set $W^\eta$, cf. Zeidler \cite{Zeidler}. 

$\rm (ii)$ We introduce the set $L=\{x\in  \Omega: \tilde u(x)\notin [0,\delta]\}$
and suppose indirectly that $m(L)>0$.   Define the function $\gamma:\mathbb
R\to \mathbb R$ by $\gamma(s)=\min(s_+,\delta)$, where
$s_+=\max(s,0).$ Now, set $w=\gamma\circ  \tilde u.$ It is clear that
$\gamma$ is a Lipschitz function and $\gamma(0)=0$. Accordingly, based on the superposition theorem of
Marcus and Mizel  \cite{MM}, one has that $w \in H_0^1(\Omega).$
Moreover,  $0\leq w(x)\leq \delta$ for a.e. $\Omega$.
Consequently, $w\in W^\eta.$

Let us introduce the sets $$L_1=\{ x\in L :
\tilde u(x)<0\}\ \ {\rm and}\ \ L_2= \{x\in L :
\tilde u(x)>\delta\}.$$ In particular, $L=L_1\cup L_2,$ and by definition, it follows that
$w(x)= \tilde u(x)$ for all $x\in \Omega\setminus L$, $w(x)=0$ for
all $x\in L_1$, and $w(x)=\delta$ for all $x\in L_2$. In addition,
one has
\begin{eqnarray*}
\mathcal T(w) -\mathcal T(\tilde u)	&=&\frac{1}{2}\left[\|w\|_{H_0^1}^2-\|\tilde
	u\|_{H_0^1}^2\right]+\frac{k}{2}\int_{\Omega}\left[w^2-\tilde
	u^2\right]-\int_{\Omega}[A(w(x))-A(\tilde u(x))]\\ &=&
	-\frac{1}{2}\int_L|\nabla \tilde u|^2+\frac{k}{2}\int_L
[w^2-\tilde u^2]-\int_L [A(w(x))-A(\tilde u(x))].
\end{eqnarray*}
On account of $k> 0$, we have  $$k\int_L [w^2- \tilde
u^2]=-k\int_{L_1}  \tilde u^2+k\int_{L_2} [\delta^2-\tilde
u^2]\leq 0.$$ Since $A(s)=0$ for all $s\leq 0,$
we have $$\int_{L_1} [A(w(x))-A(\tilde u(x))]=0.$$ By means of the Lebourg's mean value
theorem, for a.e. $x\in L_2$, there exists $\theta(x)\in [\delta,
\tilde u(x)]\subseteq [\delta,\eta]$ such that $$A(w(x))-A(
\tilde u(x))=A(\delta)-A( \tilde
u(x))=a(\theta(x))(\delta-\tilde u(x)),$$ where 
$a(\theta(x))\in \partial A(\theta(x))$. 
Due to 
{\rm (H$_A^2$)}, it turns out that $$\int_{L_2} [A(w(x))-A(\tilde u(x))]\geq
0.$$
Therefore, we obtain that $\mathcal T(w) -\mathcal T(\tilde u)\leq 0$.
 On the other hand, since $w\in W^\eta$, then  $
\mathcal T(w)\geq \mathcal T(\tilde u)=\inf_{W^\eta}\mathcal T$, thus  every term in the difference $\mathcal T(w) -\mathcal T(\tilde u)$ should be 
zero;
in particular, 
$$\int_{L_1} \tilde u^2=\int_{L_2} [\tilde
u^2-\delta^2]=0.$$ The latter relation implies in particular that $m(L)=0$, which is a contradiction,  completing the proof of (ii). 

$\rm (iii)$  Since $\tilde u(x)\in [0,\delta]$ for a.e. $x\in
\Omega$, an arbitrarily small perturbation $\tilde u+\epsilon v$ of $\tilde u$ with $0<\epsilon\ll 1$ and $v\in C_0^\infty(\Omega)$ still implies that $\mathcal T(\tilde u+ \varepsilon v)\geq \mathcal T(\tilde u)$; accordingly, $\tilde u$ is a minimum point for $\mathcal T$ in the strong topology of $H_0^1(\Omega)$, thus $0\in \partial \mathcal T(\tilde u)$, cf. Remark \ref{rem-1-0} in the Appendix. Consequently, it follows that  $\tilde u$ is a weak solution 
 of the differential inclusion $({\rm D}_A^k).$ 
\hfill $\square$\\

In the sequel, we need a truncation function of $H_0^1(\Omega)$, see also \cite{KM}.
 To construct this function, 
let $B(x_0,r)\subset \Omega$ be the $n$-dimensional ball with
radius $r>0$ and center $x_0\in \Omega$. For $s>0$, define
\begin{equation}\label{u-k-szerkesztese}
\begin{array}{l}
w_s(x) = \left \{
\begin{array}{lll}
0,& {\rm if} & \ x\in \Omega\setminus B(x_0,r);\\
s, & {\rm if} &\ x\in B(x_0,r/2);\\
\frac{2s}{r}(r-|x-x_0|), & {\rm if} & \ x\in B(x_0, r)\setminus B(x_0,r/2).
\end{array} \right. \\
\end{array}
\end{equation}
Note that  that $w_s\in
H_0^1(\Omega),$  $ \|w_s\|_{L^\infty}=s$ and
\begin{equation}\label{H01 norma, zs}
\| w_s\|_{H_0^1}^2=\int_\Omega |\nabla w_s|^2=
4r^{n-2}(1-2^{-n})\omega_n s^2\equiv C(r,n) s^2>0;
\end{equation}
hereafter $\omega_n$ stands for the volume of $B(0,1)\subset \mathbb R^n.$

\section[Nonlinearities with oscillation near the origin]{Proof of Theorems \ref{elso-tetel} and \ref{masodik-tetel}}\label{proof of theorems with osc near orig}

Before giving the proof of  Theorems \ref{elso-tetel} and \ref{masodik-tetel}, in the first part of this section we study the  differential inclusion  problem  $$\left\{
	\begin{array}{lll}
	-\triangle u(x) + ku(x) \in   \partial A(u(x)),\ \ u(x)\geq 0 \ \ &  &
	x\in \Omega,
	\\ u(x)= 0 & &
	x\in \partial \Omega,
	\end{array}\right.
	\eqno{({\rm D}_A^k)}$$
where  $k>0$ and the locally Lipschitz function $A:\mathbb R_+\to \mathbb R$ verifies
\begin{enumerate}
	\item[(H$_0^0$):]  $A(0)=0$; \item[(H$_1^0$):] 
		$-\infty<\liminf_{s\to 0^+}\frac{A(s)}{s^2}$   and
		$\limsup_{s\to 0^+}\frac{A(s)}{s^2}=+\infty;$
	\item[(H$_2^0$):]  there are two sequences $\{\delta_i\}$, $\{\eta_i\}$
	with
	$0<\eta_{i+1}<\delta_i<\eta_i$, $\lim_{i\to \infty}\eta_i=0,$ and $$\max\{\partial A(s)\}:= \max \{\xi:\xi\in \partial A(s)\}\leq 0$$
	 for every
	$s\in [\delta_i,\eta_i],$ $i\in \mathbb
	N$.
\end{enumerate}

\begin{theorem}\label{3.1-tetel} Let $k>0$ and assume hypotheses  {\rm (H$_0^0$)}, {\rm (H$_1^0$)} and  {\rm
		(H$_2^0$)} hold. Then there exists a sequence $\{
	u_i^0\}_{i}\subset H_0^1(\Omega)$ of distinct weak solutions of the differential inclusion 
	problem $({\rm D}_A^k)$ such that
	\begin{equation}\label{meg-ujabb-nulla-osc-hatarertek}
	\lim_{i\to
		\infty}\|u_i^0\|_{H_0^1}=\lim_{i\to \infty}\|
	u_i^0\|_{L^\infty}=0.
	\end{equation}
\end{theorem}

{\it Proof.} We may assume that
	$\{\delta_i\}_i$, $\{\eta_i\}_i\subset (0,1)$. For any fixed number $i\in \mathbb N$, we define the locally Lipschitz function $A_i: \mathbb R\to \mathbb R$ by
	\begin{equation}\label{trunc-f-es-g}
		A_i(s)=A(\tau_{\eta_i}(s)),
	\end{equation}
	where $A(s)=0$ for $s\leq 0$ and $\tau_\eta:\mathbb R\to \mathbb R$ denotes the truncation function 
	$\tau_\eta(s)=\min(\eta,s),$ $\eta>0.$ For further use, we introduce the energy functional $\mathcal T_i:H_0^1(\Omega)\to \mathbb R$  associated with problem $({\rm D}_{A_i}^k)$. 
	
	We notice that for $s\geq 0$, the chain rule (see Proposition \ref{chain rule1} in the Appendix) gives 
		$$\partial A_i(s)=\left\{
	\begin{array}{lll}
	\partial A(s) & {\rm if} &
	s<  \eta_i,
	\\ 
	\overline{\rm co}\{0,\partial A(\eta_i)\}& {\rm if} &
	s= \eta_i,\\
	\{0\} & {\rm if} &
	s> \eta_i.
	\end{array}\right.
	$$
	It turns out that on the compact set $[0,\eta_i]$, the upper semicontinuous set-valued map 
	$s\mapsto \partial A_i(s)$	 attains its supremum  (see Proposition \ref{prop-aaaa} in the Appendix); therefore,  there exists $M_{A_i}>0$ such that 
		$$\max |\partial A_i(s)|:=\max\{|\xi|:\xi\in \partial A_i(s)\}\leq M_{A_i}$$ 
	for every $s\geq 0$,
	i.e., (H$_{A_i}^1$) holds. The same is true for (H$_{A_i}^2$) by using 
	(H$_2^0$) on $[\delta_i,\eta_i],$ $i\in \mathbb N.$
	 
	Accordingly, the assumptions of Theorem \ref{segedtetel eredmeny oszcillacio nullaban} are verified for every
	$i\in \mathbb N$ with $[\delta_i,\eta_i],$ thus there exists $u_i^0\in W^{\eta_i}$ such that
	\begin{equation}\label{masodik-tetel-alapjan}
		u_i^0 \ \mbox{is the minimum point of the functional}\
		\mathcal T_i\ \mbox{on}\  W^{\eta_i},
	\end{equation}
	\begin{equation}\label{elso-tetel-alapjan}
		u_i^0(x)\in [0,\delta_i]\  \mbox{for a.e.}\ x\in \Omega,
	\end{equation}
	\begin{equation}\label{weak-masodik-tetel-alapjan}
		u_i^0  \ \mbox{is a  solution of}\ ({\rm D}_{A_i}^k).
	\end{equation}
		On account of relations (\ref{trunc-f-es-g}), (\ref{elso-tetel-alapjan}) and
		(\ref{weak-masodik-tetel-alapjan}), $u_i^0$ is a weak solution also for the  differential inclusion problem  $({\rm
			D}_{A}^k)$. 	

		We are going to prove that there are infinitely many distinct elements in
		the sequence $\{u_i^0\}_i$. To conclude it, we first prove that
	\begin{equation}\label{E_i negativ}
		\mathcal T_i(u_i^0)<0\ \ {\rm for\ all}\ i\in \mathbb N; {\rm\ and}
	\end{equation}
	\begin{equation}\label{E_i hatarertek}
		\lim_{i\to \infty} \mathcal T_i(u_i^0)=0.
	\end{equation}
		The left part of (H$_1^0$)
		implies the existence of some $l_0>0$ and $ \zeta\in (0,\eta_1)$ such that
	\begin{equation}\label{A also becsles vegtelen oszc}
		 A(s)\geq-l_0 s^2\ {\rm for\ all}\ s\in(0,
		\zeta).
	\end{equation}
	One can choose $L_0>0$ such that
	\begin{equation}\label{becsles L0-al}
		\frac{1}{2}C(r,n)+\left(\frac{k}{2}+l_0\right)
		m(\Omega)< L_0(r/2)^n\omega_n,
	\end{equation}
	where $r>0$ and $C(r,n)>0$ come from (\ref{H01 norma, zs}).
	Based on the right part of (H$_1^0$), one can find a sequence $\{\tilde s_i\}_i
	\subset(0, \zeta)$
	such that $\tilde s_i\leq \delta_i$ and
	\begin{equation}\label{ujabb-becsles}
		 A(\tilde s_i)>L_0 \tilde s_i^2\ \
		{\rm for\ all}\ i\in \mathbb N.
	\end{equation} Let $i\in \mathbb N$ be a fixed number and let $w_{\tilde s_i}\in
	H_0^1(\Omega)$ be the function from
	(\ref{u-k-szerkesztese}) corresponding  to the value $\tilde s_i>0.$
	Then $w_{\tilde s_i}\in
	W^{\eta_i}$,
	and due to (\ref{A also becsles vegtelen oszc}), (\ref{ujabb-becsles}) and (\ref{H01 norma, zs}) one has
	\begin{eqnarray*}
		\mathcal T_i(w_{\tilde s_i})&=&\frac{1}{2}\|w_{\tilde s_i}\|_{H_0^1}^2+
		\frac{k}{2}\int_\Omega w_{\tilde s_i}^2-\int_{\Omega}
		A_i(w_{\tilde s_i}(x))dx\\ &=&\frac{1}{2}C(r,n)\tilde s_i^2+
		\frac{k}{2}\int_\Omega w_{\tilde s_i}^2- \int_{B(x_0,r/2)}
		A({\tilde s_i})dx-\int_{B(x_0,r)\setminus B(x_0,r/2)}
		A(w_{\tilde s_i}(x))dx
		\\&\leq & \left[\frac{1}{2}C(r,n)+\frac{k}{2}m(\Omega)-
		L_0(r/2)^n\omega_n+l_0m(\Omega)\right]\tilde s_i^2.
	\end{eqnarray*}
	Accordingly, with (\ref{masodik-tetel-alapjan}) and (\ref{becsles L0-al}), we conclude that
	\begin{equation}\label{R_ies R_i hasonlintas}
		\mathcal T_i(u_i^0)=\min_{W^{\eta_i}}\mathcal T_i \leq \mathcal T_i(w_{\tilde
		s_i})<0
	\end{equation}
	which completes the proof of (\ref{E_i negativ}). 
		
	Now, we prove (\ref{E_i hatarertek}). For every $i\in \mathbb N$, by
	using the Lebourg's mean value theorem, relations  (\ref{trunc-f-es-g}) and
	(\ref{elso-tetel-alapjan}) and (H$_0^0$) , we have $$\mathcal T_i(u_i^0)\geq
	-\int_{\Omega} A_i(u_i^0(x))dx=-\int_{\Omega} A_1(u_i^0(x))dx\geq -M_{A_1}m(\Omega)\delta_i.$$ Since $\lim_{i\to
		\infty}\delta_i=0$, the latter estimate and (\ref{R_ies R_i
		hasonlintas}) provides relation (\ref{E_i hatarertek}).
	
	Based on (\ref{trunc-f-es-g}) and
	(\ref{elso-tetel-alapjan}), we have that $\mathcal
	T_i(u_i^0)=\mathcal T_1(u_i^0) \ \ {\rm for\ all}\ i\in \mathbb
	N.$ This relation with (\ref{E_i negativ}) and
	(\ref{E_i hatarertek}) means that the sequence $\{u_i^0\}_i$
	contains infinitely many distinct elements.
	
	We now  prove 
	(\ref{meg-ujabb-nulla-osc-hatarertek}). One can prove the former limit by (\ref{elso-tetel-alapjan}), i.e.
	$\|u_i^0\|_{L^\infty}\leq \delta_i$ for all $i\in \mathbb N$,
	combined with $\lim_{i\to \infty} \delta_i=0.$ For the latter
	limit, we use $k>0$, (\ref{R_ies R_i hasonlintas}), 
	(\ref{trunc-f-es-g}) and (\ref{elso-tetel-alapjan}) to get for
	all $i\in \mathbb N$ that 
	\begin{eqnarray*}
		\frac{1}{2}\|u_i^0\|_{H_0^1}^2&\leq&
		\frac{1}{2}\|u_i^0\|_{H_0^1}^2+ \frac{k}{2}\int_\Omega
		(u_i^0)^2 <\int_{\Omega} A_i(u_i^0(x))=\int_{\Omega} A_1(u_i^0(x))\leq M_{A_1} m(\Omega)\delta_i,
	\end{eqnarray*}
	which completes the proof.\hfill$\square$

{\rm \textbf{Proof of Theorem \ref{elso-tetel}}.} We split the proof into two parts.\\
	(i) {\it Case $p=1.$} Let
	$\lambda\geq 0$ with $\lambda
	\overline c<-l_0 $  and fix $\tilde\lambda_0\in \mathbb R$ such that 
	$\lambda
	\overline c<\tilde\lambda_0<-l_0.$ 
	With these choices we define   
		\begin{equation}\label{elso-valasztas}
		k:=\tilde\lambda_0-\lambda\overline c>0 \ \ {\rm and}\ A(s):=F(s)+\frac{\tilde\lambda_0
		}{2}s^2+\lambda\left(G(s)-\frac{\overline c
	}{2}s^2\right)\ {\rm  for\ every}\ s\in [0,\infty).
		\end{equation}
	It is clear that  $A(0)=0$, i.e., (H$_0^0$) is verified. 
	Since $p=1$, by $(G^0_1)$ one has  $$\underline c=\liminf_{s\to 0^+}\frac{\min\{ \partial G(s)\}}{s}\leq \limsup_{s\to 0^+}\frac{\max\{\partial G(s)\}}{s}=\overline c.$$
	In particular, for sufficiently small $\epsilon>0$ there exists $\gamma=\gamma(\epsilon)>0$ such that
	$$\max\{ \partial G(s)\}-\overline c s<\epsilon s,\ \forall s\in [0,\gamma],$$ and 
	$$\min\{ \partial G(s)\}-\underline c s>-\epsilon s,\ \forall s\in [0,\gamma].$$
	For $s\in [0,\gamma]$, Lebourg's mean value theorem and $G(0)=0$ implies that there exists $\xi_s\in \partial G(\theta_s s)$ for some $\theta_s\in [0,1]$ such that
	$G(s)-G(0)=\xi_s s$. Accordingly, for every $s\in [0,\gamma]$ we have that
	\begin{equation}\label{ket-oldalu}
	(\underline c -\epsilon) s^2\leq G(s)\leq (\overline c +\epsilon) s^2.
	\end{equation}
	By (\ref{ket-oldalu}) and $(F^0_1)$ we have that 
	$$\liminf_{s\to 0^+}\frac{A(s)}{s^2}\geq \liminf_{s\to 0^+}\frac{F(s)}{s^2}+\frac{\tilde \lambda_0-\lambda \overline c}{2}+\lambda \liminf_{s\to 0^+}\frac{G(s)}{s^2}\geq \liminf_{s\to 0^+}\frac{F(s)}{s^2}+\frac{\tilde \lambda_0-\lambda \overline c}{2}+\lambda \underline c>-\infty$$
	and 
	$$\limsup_{s\to 0^+}\frac{A(s)}{s^2}\geq \limsup_{s\to 0^+}\frac{F(s)}{s^2}+\frac{\tilde \lambda_0-\lambda \overline c}{2}+\lambda \liminf_{s\to 0^+}\frac{G(s)}{s^2}=+\infty,$$ i.e., 
 (H$_1^0$) is verified. 
	
	Since
	\begin{equation}\label{incluzio}
		\partial A(s)\subseteq \partial F(s)+\tilde \lambda_0 s +\lambda(\partial G(s)-\overline c s),
	\end{equation}
	and $\lambda\geq 0$, we have that
	\begin{equation}\label{max-becs}
		\max\{\partial A(s)\}\leq \max\{\partial F(s)+\tilde \lambda_0 s\}+\lambda\max\{\partial G(s)-\overline c s\}.
	\end{equation}
	Since $$\limsup_{s\to 0^+}\frac{\max\{\partial G(s)\}}{s}=\overline c,$$ cf. $(G^0_1)$, and 
	$$\liminf_{s\to 0^+}\frac{\max\{\partial F(s)\}}{s}=l_0<0,$$ cf. $(F^0_2)$, it turns out by (\ref{max-becs}) that
	$$\liminf_{s\to 0^+}\frac{\max\{\partial A(s)\}}{s}\leq \liminf_{s\to 0^+}\frac{\max\{\partial F(s)\}}{s}+ \tilde \lambda_0-\lambda\overline c+\lambda\limsup_{s\to 0^+}\frac{\max\{\partial G(s)\}}{s}\leq  l_0+\tilde \lambda_0<0.$$
	Therefore, one has a sequence $\{s_i\}_i\subset
	(0,1)$ converging to 0 such that $\frac{\max\{\partial A(s_i)\}}{s_i}<0$ i.e., $\max\{\partial A(s_i)\}<0$ for
	all $i\in \mathbb N.$ By using the upper semicontinuity of
	$s\mapsto \partial A(s)$,  we may choose two numbers $\delta_i,\eta_i\in (0,1)$ with $\delta_i<s_i<\eta_i$ such that $\partial A(s)\subset \partial A(s_i)+[-\epsilon_i,\epsilon_i]$ for every  $s\in
	[\delta_i,\eta_i]$, where $\epsilon_i:=-\max\{\partial A(s_i)\}/2>0.$ In particular, $\max\{\partial A(s)\}\leq 0$ for all $s\in
	[\delta_i,\eta_i]$.
	Thus, one may fix two sequences
	$\{\delta_i\}_i,\{\eta_i\}_i\subset (0,1)$ such that
	$0<\eta_{i+1}<\delta_i<s_i<\eta_i$, $\lim_{i\to \infty}\eta_i=0,$
	and $\max\{\partial A(s)\}\leq 0$ for all $s\in
	[\delta_i,\eta_i]$ and $i\in \mathbb N$. Accordingly, (H$_2^0$)
	is verified as well. Let us apply Theorem \ref{3.1-tetel} with the
	choice (\ref{elso-valasztas}), i.e., 
	there exists  a sequence $\{u_i\}_{i}\subset H_0^1(\Omega)$ of different elements such that
	 $${\left\{
		\begin{array}{lll}
		-\triangle u_i(x) + (\tilde\lambda_0-\lambda\overline c)u_i(x) \in    \partial F(u_i(x))+\tilde \lambda_0 u_i(x) +\lambda(\partial G(u_i(x))-\overline c u_i(x)) &  &
		x\in \Omega,\\
	 u_i(x)\geq 0  &  &
		x\in \Omega,
		\\
		u_i(x)= 0 & &
		x\in \partial \Omega,
		\end{array}\right.
	}$$
	where we used the inclusion (\ref{incluzio}). In particular, $u_i$ solves problem 
	$({\mathcal D}_\lambda)$, $i\in \mathbb N$, which completes the proof of (i).

	(ii) {\it Case $p>1.$} Let $\lambda\geq 0$ be arbitrary
	fixed and choose a number $\lambda_0\in (0,-l_0)$. Let 
	\begin{equation}\label{masodik-valasztas}
		k:=\lambda_0>0 \ \ {\rm and}\ A(s):=F(s)+\lambda G(s) +\lambda_0
		\frac{s^2}{2}\ {\rm  for\ every}\ s\in  [0,\infty).
	\end{equation}
	Since $F(0)=G(0)=0$, hypothesis  (H$_0^0$) clearly
	holds.	 
	By $(G^0_1)$ one has  $$\underline c=\liminf_{s\to 0^+}\frac{\min\{ \partial G(s)\}}{s^p}\leq \limsup_{s\to 0^+}\frac{\max\{\partial G(s)\}}{s^p}=\overline c.$$
	In particular, since $p>1$, then
	\begin{equation}\label{hatarertek-0}	\lim_{s\to 0^+}\frac{\min\{ \partial G(s)\}}{s}=\lim_{s\to 0^+}\frac{\max\{ \partial G(s)\}}{s}=0
	\end{equation}
	and 
	for sufficiently small $\epsilon>0$ there exists $\gamma=\gamma(\epsilon)>0$ such that
	$$\max\{ \partial G(s)\}-\overline c s^p<\epsilon s^p,\ \forall s\in [0,\gamma]$$ and 
	$$\min\{ \partial G(s)\}-\underline c s^p>-\epsilon s^p,\ \forall s\in [0,\gamma].$$
	For a fixed $s\in [0,\gamma]$, by Lebourg's mean value theorem and $G(0)=0$ we conclude again that $G(s)-G(0)=\xi_s s$. Accordingly,  for sufficiently small $\epsilon>0$ there exists $\gamma=\gamma(\epsilon)>0$ such that
	$(\underline c -\epsilon) s^{p+1}\leq G(s)\leq (\overline c +\epsilon) s^{p+1}$ for every $s\in [0,\gamma]$. Thus, since $p>1$, 
	$$\lim_{s\to 0^+}\frac{G(s)}{s^2}=\lim_{s\to 0^+}\frac{G(s)}{s^{p+1}}s^{p-1}=0.$$
	Therefore, by using (\ref{masodik-valasztas}) and ($F_1^0$), we conclude that 
	$$\liminf_{s\to 0^+}\frac{A(s)}{s^2}= \liminf_{s\to 0^+}\frac{F(s)}{s^2}+\lambda\lim_{s\to 0^+}\frac{G(s)}{s^2}+\frac{\lambda_0}{2}>-\infty,$$
	and $$\limsup_{s\to 0^+}\frac{A(s)}{s^2}=\infty,$$
	i.e., (H$_0^1$) holds. Since
	\begin{equation*}
	\partial A(s)\subseteq \partial F(s)+\lambda \partial G(s) +\lambda_0 s,
	\end{equation*}
	and $\lambda\geq 0$, we have that
	\begin{equation*}
	\max\{\partial A(s)\}\leq \max\{\partial F(s)\}+\max\{\lambda\partial G(s)+\lambda_0 s\}.
	\end{equation*}
	Since $$\limsup_{s\to 0^+}\frac{\max\{\partial G(s)\}}{s^p}=\overline c,$$
	cf. (G$_1^0$), and
	$$\liminf_{s\to 0^+}\frac{\max\{\partial F(s)\}}{s}=l_0,$$
	cf. (F$_2^0$), by relation (\ref{hatarertek-0}) it turns out that 
	$$\liminf_{s\to 0^+}\frac{\max\{\partial A(s)\}}{s}= \liminf_{s\to 0^+}\frac{\max\{\partial F(s)\}}{s} + \lambda\lim_{s\to 0^+}\frac{\max\{\partial G(s)\}}{s}+\lambda_0=  l_0+\lambda_0<0,$$ and the upper semicontinuity of $\partial A$ implies the existence of 
	two sequences
	$\{\delta_i\}_i$ and $\{\eta_i\}_i\subset (0,1)$ such that
	$0<\eta_{i+1}<\delta_i<s_i<\eta_i$, $\lim_{i\to \infty}\eta_i=0,$
	and $\max\{\partial A(s)\}\leq 0$ for all $s\in [\delta_i,\eta_i]$ and
	$i\in \mathbb N$. Therefore, 
	hypothesis (H$_2^0$) holds. Now, we can apply Theorem
	\ref{3.1-tetel}, i.e., 
	there is  a sequence $\{u_i\}_{i}\subset H_0^1(\Omega)$ of different elements such that
		$$\left\{
		\begin{array}{lll}
		-\triangle u_i(x) + \lambda_0u_i(x) \in  \partial A(u_i(x))\subseteq  \partial F(u_i(x))+\lambda\partial G(u_i(x))+ \lambda_0 u_i(x)  &  &
		x\in \Omega,\\
		u_i(x)\geq 0  &  &
		x\in \Omega,
		\\
		u_i(x)= 0 & &
		x\in \partial \Omega,
		\end{array}\right.
		$$
	which means that $u_i$ solves problem 
	$({\mathcal D}_\lambda)$, $i\in \mathbb N$. 
	This completes the proof
of Theorem \ref{elso-tetel}. \hfill $\square$ \\

\rm {\textbf{Proof of Theorem \ref{masodik-tetel}.}} The proof is done in two steps:

	(i) Let $\lambda_0 \in (0, -l_0),$ $ \lambda \geq 0$ and define  \begin{equation}\label{harmadik-valasztas}
		k:=\lambda_0>0 \ \ {\rm and}\ A^\lambda(s):=F(s)+\lambda G(s) +\lambda_0
		\frac{s^2}{2}\ {\rm  for\ every}\ s\in  [0,\infty).
	\end{equation}
One can observe that $\partial A^{\lambda}(s) \subseteq \partial F(s) + {\lambda}_0s+\lambda\partial G(s)$ for every $s\geq 0$. On  account of $(F_2^0)$, there is a sequence $\{s_i\}_i \subset (0,1)$ converging to 0 such that $$\max\{\partial A^{\lambda=0}(s_i)\} \leq \max \{\partial F(s_i)\}+\lambda_0 s_i<0.$$ Thus, due to the  upper semicontinuity of $(s,\lambda)\mapsto\partial A^\lambda(s)$, we can choose three sequences $\{\delta_i\}_i, \{\eta_i\}_i, \{\lambda_i\}_i \subset (0,1)$ such that $0<\eta_{i+1}<\delta_i<s_i<\eta_i, \lim_{i \to \infty}\eta_i = 0$, and 
	$$\max\{\partial A^\lambda(s)\} \leq 0 \ \textrm{for all}\ \lambda \in [0, \lambda_i], s \in [\delta_i, \eta_i],\ i \in \mathbb N.$$
Without any loss of generality, we may  choose
	\begin{equation}\label{delta}
		\delta_i \leq \min\{i^{-1}, 2^{-1}i^{-2}[1+ m(\Omega)(\max_{s \in[0,1]}|\partial F(s)|+\max_{s \in[0,1]}|\partial G(s)|)]^{-1}\}.
	\end{equation}
	
	For every $i \in  \mathbb N$ and $\lambda \in [0, \lambda_{i}]$, let $A_i^{\lambda}:[0,\infty) \to \mathbb R$ be defined as
	\begin{equation}\label{truncation}
		A_i^{\lambda}(s) =  A^{\lambda}(\tau_{\eta_i}(s)),
	\end{equation}
	and the energy functional $\mathcal
	T_{i,\lambda}:H_0^1(\Omega) \to \mathbb R$ associated with the differential inclusion problem$({\rm D}_{A_{i}^{\lambda}}^k)$ is given by 
		 $$\mathcal
		 T_{i,\lambda}(u)=\frac{1}{2}\|u\|^2_{H^1_0}+\frac{k}{2}\int_{\Omega}u^2dx-\int_{\Omega}A_i^{\lambda}(u(x))dx.$$
	One can easily check that for every $i \in \mathbb N$ and $\lambda \in [0, \lambda_i]$, the function $A_i^\lambda$ verifies the hypotheses of Theorem \ref{segedtetel eredmeny oszcillacio nullaban}. Accordingly, for every $i \in \mathbb N$ and $\lambda \in [0, \lambda_i]$:
	\begin{equation}\label{EiLambda min}
		 \mathcal T_{i,\lambda}\ {\rm attains\ its\ infinum\ on}\ W^{\eta_i}\ {\rm at\ some\ } u_{i,\lambda}^0 \in W^{\eta_i} 
	\end{equation}
	\begin{equation} \label{ui0 range}
	 u_{i,\lambda}^0(x) \in [0, \delta_i] {\rm\ for\ a.e.\ } x\in \Omega; 
		\end{equation}
	\begin{equation} \label{ui0 is a weak sol}
		u_{i,\lambda}^0 {\rm\ is\ a\ weak\ solution\ of\ }({\rm D}_{A_{i}^{\lambda}}^k).
	\end{equation}		
	By the choice of the function $A^\lambda$ and $k>0$, $u_{i,\lambda}^0$ is also a solution to the  differential inclusion problem $({\rm D}_{A^\lambda}^k)$, so $({\mathcal D}_\lambda)$.
		
	(ii) It is clear that for $\lambda = 0$, the set-valued map $\partial A_i^{\lambda} = \partial A_i^{0}$ verifies the hypotheses of Theorem \ref{3.1-tetel}. In particular,  $\mathcal T_i:=\mathcal T_{i,0}$ is the energy functional associated with problem $({\rm D}_{A_{i}^0}^k)$. Consequently, the elements $u_i^0:=u_{i,0}^0$ verify not only (\ref{EiLambda min})-(\ref{ui0 is a weak sol}) but also 
 \begin{equation}\label{11}
 \mathcal T_i(u_i^0) = \min_{W^{\eta_i}}\mathcal T_i\leq\mathcal T_i(w_{\tilde{s}_i})<0 {\rm \ for\ all\ } i \in \mathbb N{\rm.} 
 \end{equation}

	 Similarly to Krist\'aly and Moro\c sanu \cite{KM}, let $\{{\theta}_i\}_i$ be a sequence with negative terms such that $\lim_{i \to \infty}\theta_i = 0$. Due to (\ref{11}) we may assume that
 \begin{equation}\label{12}
			\theta_i < \mathcal T_i(u_i^0) \leq \mathcal T_i(w_{\tilde{s}_i}) <{\theta}_{i+1}. 
 \end{equation}
 Let us choose
 \begin{equation}\label{13}
	   		{\lambda}_i^{'}=\frac{\theta_{i+1}-\mathcal T_i(w_{\tilde{s}_i})}{m(\Omega)\max_{s \in[0,1]}|G(s)|+1} {\rm\ and\ }    {\lambda}_i^{''}=\frac{\mathcal T_i(u_i^0)-\theta_i}{m(\Omega)\max_{s \in[0,1]}|G(s)|+1} {\rm\ ,\ } i \in \mathbb N,   
 \end{equation}
 and for a fixed $k\in  \mathbb N$, set
  \begin{equation}
	    	{\lambda}_k^0= \min(1,\lambda_1,..., \lambda_k,{\lambda_1}^{'},...,{\lambda_k}^{'}, {\lambda_1}^{''},...,{\lambda_k}^{''}) > 0.
  \end{equation} 
  Having in our mind these choices, for every $i \in \{1,....,k\}$ and $\lambda \in[0, \lambda_k^0]$ one has
  \begin{eqnarray}\label{kelllll}
   \mathcal T_{i,\lambda}(u_{i,\lambda}^0) &\leq& \mathcal T_{i,\lambda}(w_{\tilde{s}_i})  \nonumber
   =\frac{1}{2}\|w_{\tilde{s}_i}\|_{H_0^1}^2 - \int_{\Omega}F(w_{\tilde{s}_i}(x))dx - \lambda\int_{\Omega}G(w_{\tilde{s}_i}(x))dx  \nonumber
   \\&=&\mathcal T_i(w_{\tilde{s}_i})-\lambda\int_{\Omega}G(w_{\tilde{s}_i}(x))dx \nonumber
   \\&<&\theta_{i+1},
  \end{eqnarray}   
  and due to $u_{i,\lambda}^0 \in W^{\eta i}$ and to the fact that $u_i^0$ is the minimum point of  $\mathcal T_i$ on the set $W^{\eta i}$, by (\ref{12}) we also have
  \begin{equation}\label{15}
   \mathcal T_{i,\lambda}(u_{i,\lambda}^0) = \mathcal T_i(u_{i,\lambda}^0)- \lambda\int_{\Omega}G(u_{i,\lambda}^0(x))dx \geq \mathcal T_i(u_i^0)-\lambda\int_{\Omega}G(u_{i,\lambda}^0(x))dx >\theta_i.  
  \end{equation}
	Therefore, by (\ref{kelllll}) and (\ref{15}), for every $i \in \{1,...,k\}$ and $\lambda \in [0, \lambda_k^0]$, one has
		$$ \theta_i<\mathcal T_{i,\lambda}(u_{i,\lambda}^0)<\theta_{i+1},$$ 
	 thus
	   	$$ \mathcal T_{1,\lambda}(u_{1,\lambda}^0)<{\rm...}<\mathcal T_{k,\lambda}(u_{k,\lambda}^0)<0.$$
	   	We notice that $u_i^0 \in W^{\eta_1}$ for every $i \in \{1,...,k\}$, so $\mathcal T_{i,\lambda}(u_{i,\lambda}^0) = \mathcal T_{1,\lambda}(u_{i,\lambda}^0)$ because of (\ref{truncation}). Therefore, we conclude that for every $\lambda \in [0,\lambda_k^0]$,
	   	$$\mathcal T_{1,\lambda}(u_{1,\lambda}^0)<{\rm...}<\mathcal T_{1,\lambda}(u_{k,\lambda}^0)<0=\mathcal T_{1,\lambda}(0).$$
		Based on these inequalities, it turns out that the elements $u_{1,\lambda}^0,...,u_{k,\lambda}^0$ are distinct and non-trivial whenever $\lambda \in [0, \lambda_k^0]$.
		
	 Now, we are going to prove the estimate  (\ref{becsles-i-vel}).
		We have for every $i \in \{1,...,k\}$ and $\lambda \in [0, \lambda_k^0]$:
		$$\mathcal T_{1,\lambda}(u_{i,\lambda}^0) = \mathcal T_{i,\lambda}(u_{i,\lambda}^0) < \theta_{i+1}<0.$$
		By Lebourg's mean value theorem and (\ref{delta}), we have for every $i \in \{1,...,k\}$ and $\lambda \in [0, \lambda_k^0]$  that
		\begin{eqnarray}\label{16}
			\frac{1}{2}\|u_{i,\lambda}^0\|_{H_0^1}^2&<&\int_{\Omega}F(u_{i,\lambda}^0(x))dx + \lambda\int_{\Omega}G(u_{i,\lambda}^0(x))dx\nonumber
			\\&\leq&  m(\Omega)\delta_i[\max_{s \in[0,1]}|\partial F(s)|+\max_{s \in[0,1]}|\partial G(s)|]\nonumber\\
			&\leq & \frac{1}{2i^2}.\nonumber
		\end{eqnarray}
	 This completes the proof of Theorem \ref{masodik-tetel}. \hfill $\square$

 \section[Nonlinearities with oscillation at the infinity]{Proof of Theorems \ref{harmadik tetel} and \ref{negyedik tetel}}\label{proof of theorems with osc at infinity}
 
 We consider again  the differential inclusion problem 
 $$\left\{
 \begin{array}{lll}
 -\triangle u(x) + ku(x) \in   \partial A(u(x)),\ \ u(x)\geq 0 \ \ &  &
 x\in \Omega,
 \\ u(x)= 0 & &
 x\in \partial \Omega,
 \end{array}\right.
 \eqno{({\rm D}_A^k)},$$
  where $k > 0$ and the locally Lipschitz function $A:\mathbb R_+\to \mathbb R$ verifies
\begin{enumerate}
	\item[(H$_0^{\infty}$):]  $A(0)=0$;
	\item[(H$_1^{\infty}$):] 
	$-\infty<\liminf_{s\to {\infty}}\frac{A(s)}{s^2}$ and
	$\limsup_{s\to {\infty}}\frac{A(s)}{s^2}=+\infty;$
	\item[(H$_2^{\infty}$):]  there are two sequences $\{\delta_i\}$, $\{\eta_i\}$
	with
	$0<\delta_i<\eta_i<\delta_{i+1}$, $\lim_{i\to \infty}\delta_i=\infty,$ and $$\max\{\partial A(s)\}:= \max \{\xi:\xi\in \partial A(s)\}\leq 0$$
	for every
	$s\in [\delta_i,\eta_i],$ $i\in \mathbb
	N$.
\end{enumerate}

The counterpart of Theorem \ref{3.1-tetel} reads as follows. 

\begin{theorem}\label{segedtetel vegtelen oszc}
	 Let $k>0$ and assume the hypotheses  {\rm (H$_0^{\infty}$)}, {\rm (H$_1^{\infty}$)} and  {\rm (H$_2^{\infty}$)} hold. Then the differential inclusion problem $({\rm D}_A^k)$ admits a sequence $\{u_i^{\infty}\}_{i}\subset H_0^1(\Omega)$ of distinct
	weak solutions such that
	\begin{equation}\label{nulla-osc-hatarertek-inf}
	\lim_{i\to \infty}\|u_i^{\infty}\|_{L^\infty}=\infty.
	\end{equation}
\end{theorem}

 \textit{ Proof}. The proof is similar to the one performed in  Theorem \ref{3.1-tetel}; we shall show the differences only. We associate the energy functional $\mathcal T_i:H_0^1(\Omega)\to \mathbb R$ with problem $({\rm D}_{A_i}^k)$, where $A_i: \mathbb R\to \mathbb R$ is given by
\begin{equation}\label{trunc-f-es-g-inf}
A_i(s)=A(\tau_{\eta_i}(s)),
\end{equation}
with $A(s)=0$ for $s\leq 0$. One can show  that 
there exists  $M_{A_i}>0$ such that
$$\max|\partial A_i(s)|:= \max \{|\xi|:\xi\in \partial A_i(s)\}\leq M_{A_i}$$ for all $s \geq 0$, i.e, hypothesis (H$_{A_i}^1)$ holds. Moreover, (H$_{A_i}^2)$ follows by (H$_2^{\infty})$. Thus  Theorem \ref{3.1-tetel} can be applied for all $i \in \mathbb N$, i.e., 
 we have an element $u_i^{\infty} \in W^{\eta_i}$ such that
	\begin{equation}\label{segedtetel hasznalata vegtelen oszc 1}
	u_i^{\infty} \ \mbox{is the minimum point of the functional}\
	\mathcal T_i\ \mbox{on}\  W^{\eta_i},
	\end{equation}
	\begin{equation}\label{segedtetel hasznalata vegtelen oszc 2}
	u_i^{\infty}(x)\in [0,\delta_i]\  \mbox{for a.e.}\ x\in \Omega,
	\end{equation}
	\begin{equation}\label{segedtetel hasznalata vegtelen oszc 3}
	 \textnormal{ $u_i^{\infty}$ is a weak solution of } ({\rm D}_{A_{i}}^k).
	\end{equation}
By (\ref{trunc-f-es-g-inf}),  $u_i^{\infty}$ turns to be a weak solution also for  differential inclusion problem $({\rm D}_{A}^k)$. 
	
	We shall prove that there are infinitely many distinct elements in the sequence $\{u_i^\infty\}_i$ by showing that
	 \begin{equation}\label{E limit -infinite}
		\lim_{i \to \infty} \mathcal T_i(u_i^{\infty})= -\infty.
	 \end{equation}
	By the left part of (H$_1^{\infty}$) we can find $l_{\infty}^A>0$ and $\zeta>0$ such that
	\begin{equation}\label{elso A becsles inf}
		A(s) \geq -l_{\infty}^A  \textnormal{ for all $s>\zeta$.} 
	\end{equation}
	Let us choose  $L_{\infty}^A > 0$ large enough such that
	\begin{equation}\label{masodik A becsles inf}
		\frac{1}{2}C(r,n)+\left(\frac{k}{2}+l_{\infty}^A\right)
		m(\Omega)< L_{\infty}^A(r/2)^n\omega_n.
	\end{equation}
	 On account of the right part of (H$_1^{\infty}$), one can fix a sequence $\{\tilde{s_i}\}_i \subset (0, \infty)$ such that $\lim_{i \to {\infty}}\tilde{s_i} = \infty$ and 
	\begin{equation}\label{asasa}
	A(\tilde{s_i}) > L_{\infty}^A {\tilde{s_i}}^2 \textnormal{ for every $i \in \mathbb N$.}
	\end{equation}
 We know from (H$_2^{\infty}$) that $\lim_{i \to {\infty}}{\delta_i} = \infty$, therefore one has a subsequence $\{\delta_{m_i}\}_i$ of $\{\delta_i\}_i$ such that $\tilde{s_i} \leq \delta_{m_i}$ for all $i \in \mathbb N$. Let $i \in \mathbb N$, and recall ${w_s}_i \in H_0^1(\Omega)$ from (\ref{u-k-szerkesztese}) with $s_i:={\tilde{s}}_i >0$. Then ${w_{\tilde{s}}}_i \in W^{{{\eta}_m}_i}$ and according to (\ref{H01 norma, zs}), (\ref{elso A becsles inf}) and (\ref{asasa}) we have
\begin{eqnarray*}
 	{\mathcal T_m}_i(w_{\tilde s_i})&=&\frac{1}{2}\|w_{\tilde s_i}\|_{H_0^1}^2+
	\frac{k}{2}\int_\Omega w_{\tilde s_i}^2-\int_{\Omega}
	A_{m_i}(w_{\tilde s_i}(x))dx \\
	&=&\frac{1}{2}C(r,n)\tilde s_i^2 +\frac{k}{2}\int_\Omega w_{\tilde s_i}^2- \int_{B(x_0,r/2)}
	A({\tilde s_i})dx \\
 	&&-\int_{(B(x_0,r)\setminus B(x_0,r/2))\cap \{w_{{\tilde{s}}_i} > \zeta\}}A(w_{\tilde s_i}(x))dx\\
 	 &&- \int_{(B(x_0,r)\setminus B(x_0,r/2))\cap \{w_{{\tilde{s}}_i} \leq \zeta\}} A(w_{\tilde s_i}(x))dx \\
 	&\leq&  \left[\frac{1}{2}C(r,n)+\frac{k}{2}m(\Omega)-
 	L_{\infty}^A(r/2)^n\omega_n+l_{\infty}^Am(\Omega)\right]\tilde s_i^2 + \tilde M_A m({\Omega})\zeta,
\end{eqnarray*}
where $\tilde M_A=\max\{|A(s)|:s\in [0,\zeta]\}$ does not depend on $i\in \mathbb N.$
This estimate combined by (\ref{masodik A becsles inf}) and $\lim_{i \to {\infty}}{\tilde{s}}_i=\infty$ yields that
\begin{equation}\label{limit-qqqq}
	\lim_{i \to \infty} {\mathcal T_{m_i}} (w_{{\tilde{s}}_i}) = -\infty. 
\end{equation}
By equation (\ref{segedtetel hasznalata vegtelen oszc 1}), one has 
\begin{equation}
	{\mathcal T_{m_i}} (u_{m_i}^{\infty}) = \min_{W^{{\eta_m}_i}}{\mathcal T_{m_i}} \leq {\mathcal T_{m_i}}(w_{{\tilde{s}}_i}). 
\end{equation}
It follows by (\ref{limit-qqqq}) that 
$
\lim_{i \to \infty} {\mathcal T_{m_i}} (u_{m_i}^{\infty}) = -\infty.
$

We notice that the sequence $\{\mathcal T_i(u_i^{\infty})\}_i$ is non-increasing. Indeed, let $i<k$; due to  (\ref{trunc-f-es-g-inf})  one has that
\begin{equation}\label{E becsles}
	{\mathcal T}_i (u_i^{\infty}) = \min_{W^{\eta_i}}{\mathcal T_i} = \min_{W^{\eta_i}}{\mathcal T_k} \geq \min_{W^{\eta_k}}{\mathcal T_k} = {\mathcal T_k}(u_k^{\infty}),  
\end{equation}
which completes the proof of (\ref{E limit -infinite}).

The proof of  (\ref{nulla-osc-hatarertek-inf}) goes in a similar way as in \cite{KM}. 
	\hfill $\square$\\

\rm {\textbf{Proof of Theorem \ref{harmadik tetel}.}} We split the proof into two parts. \\
(i) {\it Case $p=1.$} Let
$\lambda \geq 0$ with $\lambda
\overline c<-l_\infty $  and fix $\tilde\lambda_\infty \in \mathbb R$ such that 
$\lambda
\overline c<\tilde\lambda_{\infty}<-l_{\infty}.$ 
With these choices, we define  
\begin{equation}\label{negyedik-valasztas}
k:=\tilde\lambda_{\infty}-\lambda\overline c>0 \ \ {\rm and}\ A(s):=F(s)+\frac{\tilde\lambda_{\infty}
}{2}s^2+\lambda\left(G(s)-\frac{\overline c
}{2}s^2\right)\ {\rm  for\ every}\ s\in [0,\infty).
\end{equation}
It is clear that  $A(0)=0$, i.e., (H$_0^{\infty}$) is verified. A similar argument for the $p$-order perturbation  $\partial G$ as before shows that 
$$\liminf_{s\to \infty}\frac{A(s)}{s^2}\geq \liminf_{s\to \infty}\frac{F(s)}{s^2}+\frac{\tilde \lambda_\infty-\lambda \overline c}{2}+\lambda \liminf_{s\to \infty}\frac{G(s)}{s^2}\geq \liminf_{s\to \infty}\frac{F(s)}{s^2}+\frac{\tilde \lambda_\infty-\lambda \overline c}{2}+\lambda \underline c>-\infty,$$
and 
$$\limsup_{s\to \infty}\frac{A(s)}{s^2}\geq \limsup_{s\to \infty}\frac{F(s)}{s^2}+\frac{\tilde \lambda_\infty-\lambda \overline c}{2}+\lambda \liminf_{s\to \infty}\frac{G(s)}{s^2}=+\infty,$$
 i.e., 
(H$_1^{\infty}$) is verified.

Since
\begin{equation}\label{incluzio p = 1, inf}
\partial A(s)\subseteq \partial F(s)+\tilde \lambda_{\infty} s +\lambda(\partial G(s)-\overline c s),\ \ s\geq 0,
\end{equation}
 it turns out that
$$\liminf_{s\to \infty}\frac{\max\{\partial A(s)\}}{s}\leq \liminf_{s\to \infty}\frac{\max\{\partial F(s)\}}{s}+ \tilde \lambda_\infty-\lambda\overline c+\lambda\limsup_{s\to \infty}\frac{\max\{\partial G(s)\}}{s}=  l_\infty+\tilde \lambda_\infty<0.$$

 By using the upper semicontinuity of
$s\mapsto \partial A(s)$,   one may fix two sequences
$\{\delta_i\}_i,\{\eta_i\}_i\subset (0,\infty)$ such that
$0<\delta_i<s_i<\eta_i<\delta_{i+1}$, $\lim_{i\to \infty}\delta_i=\infty,$
and $\max\{\partial A(s)\}\leq 0$ for all $s\in
[\delta_i,\eta_i]$ and $i\in \mathbb N$. Thus, (H$_2^{\infty}$)
is verified as well. By applying the inclusion (\ref{incluzio p = 1, inf}) and Theorem \ref{3.1-tetel} with the
choice (\ref{negyedik-valasztas}),  
there exists  a sequence $\{u_i\}_{i}\subset H_0^1(\Omega)$ of different elements such that
$${\left\{
\begin{array}{lll}
-\triangle u_i(x) + (\tilde\lambda_{\infty}-\lambda\overline c)u_i(x) \in \partial F(u_i(x))+\tilde \lambda_{\infty} u_i(x) +\lambda(\partial G(u_i(x))-\overline c u_i(x)) &  &
x\in \Omega,\\
u_i(x)\geq 0  &  &
x\in \Omega,
\\
u_i(x)= 0 & &
x\in \partial \Omega,
\end{array}\right.
}$$
i.e., $u_i$ solves problem 
$({\mathcal D}_\lambda)$, $i\in \mathbb N$. 

	(ii) {\it Case $p<1.$} Let $\lambda\geq 0$ be arbitrary
fixed and choose a number $\lambda_{\infty}\in (0,-l_{\infty})$. Let 
\begin{equation}\label{otodik-valasztas}
k:=\lambda_{\infty}>0 \ \ {\rm and}\ A(s):=F(s)+\lambda G(s) +\lambda_{\infty}
\frac{s^2}{2}\ {\rm  for\ every}\ s\in  [0,\infty).
\end{equation}
Since $F(0)=G(0)=0$, hypothesis  (H$_0^{\infty}$) clearly
holds. Moreover, by  $(G_1^{\infty})$,  for sufficiently small $\epsilon>0$ there exists $s_0>0$, such that
$(\underline c -\epsilon) s^{p+1}\leq G(s)\leq (\overline c +\epsilon) s^{p+1}$ for every $s>s_0$. Thus, since $p<1$, 
$$\lim_{s\to {\infty}}\frac{G(s)}{s^2}=\lim_{s\to {\infty}}\frac{G(s)}{s^{p+1}}s^{p-1}=0.$$
Accordingly, by using (\ref{otodik-valasztas}) we obtain that 
hypothesis (H$_1^{\infty}$) holds. A similar argument as above implies that 
$$\liminf_{s\to {\infty}}\frac{\max\{\partial A(s)\}}{s}\leq l_0+\lambda_{\infty}<0,$$ and the upper semicontinuity of $\partial A$ implies the existence of 
two sequences
$\{\delta_i\}_i$ and $\{\eta_i\}_i\subset (0,1)$ such that
$0<\delta_i<s_i<\eta_i<\delta_{i+1}$, $\lim_{i\to \infty}\delta_i=\infty,$
and $\max\{\partial A(s)\}\leq 0$ for all $s\in [\delta_i,\eta_i]$ and
$i\in \mathbb N$. Therefore, 
hypothesis (H$_2^{\infty}$) holds. Now, we can apply Theorem
\ref{3.1-tetel}, i.e., 
there is  a sequence $\{u_i\}_{i}\subset H_0^1(\Omega)$ of different elements such that
$${\left\{
\begin{array}{lll}
-\triangle u_i(x) + \lambda_{\infty}u_i(x) \in  \partial A(u_i(x))\subseteq  \partial F(u_i(x))+\lambda\partial G(u_i(x))+ \lambda_{\infty} u_i(x)  &  &
x\in \Omega,\\
u_i(x)\geq 0  &  &
x\in \Omega,
\\
u_i(x)= 0 & &
x\in \partial \Omega,
\end{array}\right.
}$$
which means that $u_i$ solves problem 
$({\mathcal D}_\lambda)$, $i\in \mathbb N$, which  completes the proof. \hfill $\square$ 
\\
\\

{\rm \textbf{Proof of Theorem \ref{negyedik tetel}.}} The proof is done in two steps:

(i) Let $\lambda_{\infty} \in (0, -l_\infty), \lambda \geq 0$ and define  
\begin{equation}\label{hatodik-valasztas}
k:=\lambda_{\infty}>0 \ \ {\rm and}\ A^\lambda(s):=F(s)+\lambda G(s) +\lambda_{\infty}
\frac{s^2}{2}\ {\rm  for\ every}\ s\in  [0,\infty).
\end{equation}
One has clearly that  $\partial A^{\lambda}(s) \subseteq \partial F(s) + {\lambda}_{\infty}s+\lambda\partial G(s)$ for every $s\in \mathbb R$. On account of $(F_2^{\infty})$, there is a sequence $\{s_i\}_i \subset (0,\infty)$ converging to $\infty$ such that $$\max\{\partial A^{\lambda=0}(s_i)\} \leq \max \{\partial F(s_i)\}+\lambda_{\infty} s_i<0.$$ By the  upper semicontinuity of $(s,\lambda)\mapsto \partial A^\lambda(s)$, we can choose the sequences $\{\delta_i\}_i, \{\eta_i\}_i, \{\lambda_i\}_i \subset (0,\infty)$ such that $0<\delta_i<s_i<\eta_i<\delta_{i+1}, \lim_{i \to \infty}\delta_i = \infty$, and 
$$\max\{\partial A^\lambda(s)\} \leq 0$$ for all $\lambda \in [0, \lambda_i], s \in [\delta_i, \eta_i]$ and  $i \in \mathbb N.$

For every $i \in  \mathbb N$ and $\lambda \in [0, \lambda_{i}]$, let $A_i^{\lambda}:[0,\infty) \to \mathbb R$ be defined by  
\begin{equation}\label{truncation inf 1}
A_i^{\lambda}(s) =  A^{\lambda}(\tau_{\eta_i}(s)),
\end{equation}
and accordingly, the energy functional $\mathcal
T_{i,\lambda}:H_0^1(\Omega) \to \mathbb R$ associated with the differential inclusion problem$({\rm D}_{A_{i}^{\lambda}}^k)$ is 
$$\mathcal
T_{i,\lambda}(u)=\frac{1}{2}\|u\|^2_{H^1_0}+\frac{k}{2}\int_{\Omega}u^2dx-\int_{\Omega}A_i^{\lambda}(u(x))dx.$$
Then for every $i \in \mathbb N$ and $\lambda \in [0, \lambda_i]$, the function $A_i^\lambda$ clearly verifies the hypotheses of Theorem \ref{segedtetel eredmeny oszcillacio nullaban}. Accordingly, for every $i \in \mathbb N$ and $\lambda \in [0, \lambda_i]$ there exists
\begin{equation}\label{EiLambda min inf 1}
\mathcal T_{i,\lambda}\ {\rm attains\ its\ infimum\ at\ some\ } \tilde{u}_{i,\lambda}^{\infty} \in W^{\eta_i} 
\end{equation}
\begin{equation} \label{ui0 range inf 1}
\tilde{u}_{i,\lambda}^{\infty} \in [0, \delta_i] {\rm\ for\ a.e.\ } x\in \Omega; 
\end{equation}
\begin{equation} \label{ui0 is a weak sol inf 1}
\tilde{u}_{i,\lambda}^{\infty}(x) {\rm\ is\ a\ weak\ solution\ of\ }({\rm D}_{A_{i}^{\lambda}}^k).
\end{equation}		
Due to (\ref{truncation inf 1}),  $\tilde{u}_{i,\lambda}^{\infty}$ is not only a solution to $({\rm D}_{A_{i}^{\lambda}}^k)$ but also to the differential inclusion problem $({\rm D}_{A^\lambda}^k)$, so $({\mathcal D}_\lambda).$

(ii) For $\lambda = 0$, the function $\partial A_i^{\lambda} = \partial A_i^{0}$ verifies the hypotheses of Theorem \ref{3.1-tetel}. Moreover, $\mathcal T_i:=\mathcal T_{i,0}$ is the energy functional associated with problem $({\rm D}_{A_{i}^0}^k)$. Consequently, the elements $u_i^{\infty}:=u_{i,0}^{\infty}$ verify not only (\ref{EiLambda min inf 1})-(\ref{ui0 is a weak sol inf 1}) but also 
\begin{equation}\label{E mi min}
\mathcal T_{m_i}(u_{m_i}^{\infty}) = \min_{W^{\eta_{m_i}}}(\mathcal T_{m_i})\leq  \mathcal T_{m_i}(w_{\tilde{s}_i}) {\rm \ for\ all\ } i \in \mathbb N{\rm,} 
\end{equation}
where the subsequence $\{u_{m_i}^{\infty}\}_i$ of $\{u_i^{\infty}\}_i$  and $w_{\tilde{s}_i} \in W^{\eta_i}$ appear in the proof of Theorem \ref{segedtetel vegtelen oszc}. 

  Similarly to Krist\'aly and Moro\c sanu \cite{KM}, let $\{{\theta}_i\}_i$ be a sequence with negative terms such that $\lim_{i \to \infty}\theta_i = -\infty$. On account of (\ref{E mi min}) we may assume that
\begin{equation}\label{theta series}
\theta_{i+1} < \mathcal T_{m_i}(u_{m_i}^{\infty}) \leq \mathcal T_{m_i}(w_{\tilde{s}_i}) <{\theta}_i. 
\end{equation}
Let 
\begin{equation}\label{lambda estimation inf}
{\lambda}_i^{'}=\frac{\theta_i-\mathcal T_{m_i}(w_{\tilde{s}_i})}{m(\Omega)\max_{s \in[0,1]}|G(s)|+1} {\rm\ and\ }    {\lambda}_i^{''}=\frac{\mathcal T_{m_i}(u_{m_i}^{\infty})-\theta_{i+1}}{m(\Omega)\max_{s \in[0,1]}|G(s)|+1} {\rm\ ,\ } i \in \mathbb N,   
\end{equation}
and for a fixed $k\in  \mathbb N$, we set
\begin{equation}\label{lambda _k}
{\lambda}_k^{\infty}= \min(1,\lambda_1,..., \lambda_k,{\lambda_1}^{'},...,{\lambda_k}^{'}, {\lambda_1}^{''},...,{\lambda_k}^{''}) > 0.
\end{equation} 
Then, for every $i \in \{1,....,k\}$ and $\lambda \in[0, \lambda_k^{\infty}]$, due to (\ref{theta series}) we have that
\begin{eqnarray}
\mathcal T_{m_i,\lambda}(\tilde{u}_{m_i,\lambda}^{\infty}) &\leq& \mathcal T_{m_i,\lambda}(w_{\tilde{s}_i})  \nonumber
=\frac{1}{2}\|w_{\tilde{s}_i}\|_{H_0^1}^2 - \int_{\Omega}F(w_{\tilde{s}_i}(x))dx - \lambda\int_{\Omega}G(w_{\tilde{s}_i}(x))dx  \nonumber
\\&=&\mathcal T_{m_i}(w_{\tilde{s}_i})-\lambda\int_{\Omega}G(w_{\tilde{s}_i}(x))dx \nonumber
\\&<&\theta_i.
\end{eqnarray}   
Similarly, since $\tilde{u}_{m_i,\lambda}^{\infty} \in W^{\eta_{m_i}}$ and $u_{m_i}^{\infty}$ is the minimum point of  $\mathcal T_i$ on the set $W^{\eta_{m_i}}$, on account of (\ref{theta series}) we have
\begin{equation}
\mathcal T_{m_i,\lambda}(\tilde{u}_{m_i,\lambda}^{\infty}) = \mathcal T_{m_i}(\tilde{u}_{m_i,\lambda}^{\infty})- \lambda\int_{\Omega}G(\tilde{u}_{m_i,\lambda}^{\infty})dx \geq \mathcal T_{m_i}(u_{m_i}^{\infty})-\lambda\int_{\Omega}G(\tilde{u}_{m_i,\lambda}^{\infty})dx >\theta_{i+1}.  
\end{equation} 
Therefore, for every $i \in \{1,...,k\}$ and $\lambda \in [0, \lambda_k^{\infty}]$,
\begin{equation}\label{ThetaSeries}
\theta_{i+1}<\mathcal T_{m_i,\lambda}(\tilde{u}_{m_i,\lambda}^{\infty})<\theta_i <0, 
\end{equation}  
 thus
\begin{equation}\label{Ordered E}
\mathcal T_{m_k,\lambda}(\tilde{u}_{m_k,\lambda}^{\infty})<{\rm...}<\mathcal T_{m_1,\lambda}(\tilde{u}_{m_1,\lambda}^{\infty})<0.
\end{equation}
Because of (\ref{truncation inf 1}), we notice that  $\tilde{u}_{m_i,\lambda}^{\infty} \in W^{\eta_{m_k}}$ for every $i \in \{1,...,k\}$, thus $\mathcal T_{m_i,\lambda}(\tilde{u}_{m_i,\lambda}^{\infty}) = \mathcal T_{m_k,\lambda}(\tilde{u}_{i,\lambda}^{\infty})$. Therefore,  for every $\lambda \in [0,\lambda_k^{\infty}]$,
$$\mathcal T_{m_k,\lambda}(\tilde{u}_{m_k,\lambda}^{\infty})<{\rm...}<\mathcal T_{m_k,\lambda}(\tilde{u}_{m_1,\lambda}^{\infty})<0=\mathcal T_{m_k,\lambda}(0),$$
i.e, the elements $\tilde u_{m_1,\lambda}^{\infty},...,\tilde u_{m_k,\lambda}^{\infty}$ are distinct and non-trivial whenever $\lambda \in [0, \lambda_k^{\infty}]$.
The estimate (\ref{becsles-i-vel-infty}) follows in a similar manner as in  \cite{KM}. \hfill $\square$

\section{Concluding remarks} \label{conc-section}

\begin{itemize}
	\item[1.] Suitable modification of our arguments provide multiplicity results for the differential inclusion problem 
	\[ \   \left\{ \begin{array}{lll}
	-\Delta u(x)+u(x)\in \partial F(u(x))+\lambda \partial G(u(x))& {\rm in} &   \mathbb R^n; \\
	u\geq 0, &\mbox{in} &   \mathbb R^n,\\
	\end{array}\right. \eqno{(\tilde{\mathcal D}_\lambda)}\]
	where  $\partial F$ and $\partial G$ behave in a similar manner as before. The main difficulty in the investigation of $(\tilde{\mathcal D}_\lambda)$ is the lack of compact  embedding of the Sobolev space $H^1(\mathbb R^n)$ into the Lebesgue spaces $L^q(\mathbb R^n)$, $n\geq2$, $q\in [2,2^*)$. However, by using Strauss-type estimates and Lions-type embedding results for radially symmetric functions of $H^1(\mathbb R^n)$ (see e.g. Willem \cite{Willem}), the principle of symmetric criticality for non-smooth functionals (see Kobayashi and \^Otani \cite{KO} and Squassina \cite{Squassina}) provides the expected results.  A related result in the smooth setting can be found in Krist\'aly \cite{K-1}. 
	               
	\item[2.] Assume that $\partial F$ oscillates at a point $l\in [0,+\infty]$ and  $\partial G$ has a $p$-order growth at $l.$ We are wondering if our results, valid for $l=0$ and $l=+\infty$, can be extended to any $l\in (0,\infty)$, even in the smooth framework. 
\end{itemize}

	\section{Appendix: Locally Lipschitz functions}\label{appendix}

In this part we collect those notions and properties of locally Lipschitz functions which are used in the proofs; for details, see Clarke \cite{Clarke} and Chang \cite{Chang}. 
		Let $(X, \|\cdot\|)$ be a real Banach space and  $U\subset X$ be an open set; we denote by $\langle \cdot, \cdot \rangle$ the duality mapping between
	$X^{\star}$ and $X$.

	\begin{definition} {\rm (see \cite{Clarke})}
		A function $f:X\to\mathbb{R}$ is {{\rm locally Lipschitz}} if,
		for every $x\in X$, there exist a neighborhood $U$ of $x$ and a
		constant $L>0$ such that
		$$ |f(x_1)-f(x_2)|\leq L\|x_1-x_2\|\ \ for\ all \ \ x_1,x_2\in U.$$
	\end{definition}

	\begin{definition} {\rm (see \cite{Clarke})}
		Let  $f$ be a locally Lipschitz function near the point $x$ and let $v$ be any arbitrary vector in $X$. The {{\rm generalized directional derivative in the sense of Clarke}} of $f$ at the
		point $x \in X$ in  the direction $v \in X$ is
		$$f^\circ(x;v)=\limsup_{z\to x, \ \tau\to 0^+}\frac{f(z+\tau v)-f(z)}{\tau}.$$
		The {{\rm generalized gradient}} of $f$ at $x\in X$ is the set
		$$\partial f(x)=\{x^\star \in X^\star : \ \langle x^\star,v\rangle\leq
		f^\circ(x;v) \ \mbox{\rm for all} \ v\in X\}.$$
	\end{definition}
	
	For all $x\in X$, the functional $f^\circ(x,\cdot)$ is finite
	and positively homogeneous. Moreover, we have the following properties.
	
	\begin{proposition}\label{h} {\rm (see \cite{Clarke})} \label{prop-aaaa}
		Let $X$ be a real Banach space, $U\subset X$ an open subset and
		$f,g:U\rightarrow\mathbb{R}$ be locally Lipschitz functions.
		The following properties hold: 
		\begin{itemize}
			\item[$(a)$] For every $x\in U$, $\partial f(x)$ is a nonempty, convex
			and weakly$^{\star}$-compact subset of $X^{\star}$ which is bounded
			by the Lipschitz constant $L > 0$ of $f$ near $x;$
			\item[$(b)$] $f^\circ(x;v)=\max\{\langle \xi, v\rangle : \ \xi\in\partial
			f(x)\}$ for all $v \in X$;
			\item[$(c)$] $(f+g)^\circ(x;v)\leq  f^\circ(x;v)+g^\circ(x;v)$ for all $x\in U,$ $v\in X;$
			\item[$(d)$]	$ \partial (f+g)(u)\subset \partial f(u)+\partial g(u)$ for all\ $u\in U$;
			\item[$(e)$] $(-f)^\circ(x;v)=f^\circ(x;-v)$ for all $x \in U$;
			\item [$(f)$] The function $(x,v)\mapsto f^\circ(x;v)$ is upper
			semicontinuous;
			\item[$(g)$] The set-valued map
			$\partial f:U\rightarrow 2^{X^{\star}}$ is
			weakly$^{\star}$-closed, that is, if  $\{x_i\}\subset U$ and
			$\{w_i\}\subset X^{\ast}$ are sequences such that
			$x_i\rightarrow x$ strongly in $ X$ and $w_i\in
			\partial f(x_i)$ with $w_i \rightharpoonup z$ weakly$^{\star}$ in
			$X^{\ast}$, then $z\in \partial f(x).$
			In particular, if $X$ is finite dimensional, then $\partial f$ is
			upper semicontinuous, i.e., for every $\epsilon>0$ there exists $\gamma>0$ such that
			$\partial f(x')\subseteq \partial f(x)+ B_{X^*}(0,\epsilon),\ \forall x'\in B_X(x,\gamma);$
		\end{itemize}	
	\end{proposition}
	
	
	\begin{proposition}\label{lambda} {\rm (see \cite{Chang})} 
		The number $\lambda_f(u)= \displaystyle \inf_{w \in \partial f(u)}
		||w||_{X^{\star}}$  is well defined and  $$\displaystyle \liminf_{u \rightarrow u_{0}} \lambda_f(u) \geq
		\lambda_f(u_{0}).$$
	\end{proposition}
	
	\begin{definition}\label{kritpont-chang} {\rm (see \cite{Chang})}
		Let $f: X\rightarrow \mathbb{R}$ be a locally Lipschitz function.
		We say that  $u \in X$ is a {\rm critical point} $($in the sense of Chang$)$ of $f$, if $\lambda_f(u)=0$, i.e.,  $0 \in
		\partial f(u)$.
	\end{definition}

	\begin{remark} {\rm (see \cite{Clarke})}\rm \label{rem-1-0}
		(a) $u\in X$ is a critical point of $f$ if $f^\circ(u;v)\geq 0 \
		\mbox{for all} \ v\in X.$
		
		(b) 	If  $x\in U$ is a local minimum or maximum of the locally
		Lipschitz function $f: X \rightarrow \mathbb{R}$ on an open set of
		a Banach space, then $x$ is a critical point of $f.$
	\end{remark}

	\begin{proposition}\label{lebourg kozepertek tetel}  {\rm (see \cite{Clarke})}		$(${\rm Lebourg's mean value theorem}$)$
		Let X be a Banach space, $x,y \in X$ and $f:X \to \mathbb R$ be Lipschitz on an open set containing the line segment $[x,y]$. Then there is a point $a \in (x,y)$ such that 
		$$f(y) - f(x) \in \langle\partial f(a), y-x\rangle. $$
	\end{proposition}

	\begin{proposition}\label{chain rule1} {\rm (see \cite{Clarke})} {\rm (Chain Rule)}
		Let X be Banach space, let us consider the   composite function $f = g \circ h$ where $h: X \to \mathbb R^n$ and $g: \mathbb R^n \to \mathbb R$ are given functions. Let denote $h_i$, $i \in \{1,...,n\}$ be the component functions of $h$. We assume $h_i$ is locally Lipschitz near $x$  and $g$ is too near $h(x)$. Then $f$ is locally Lipschitz near $x$ as well. Let us denote by $\alpha_i$ the elements of $\partial g$, and let  $\alpha = (\alpha_1, ..., \alpha_n)$; then 
		$$\partial f(x)\subset \overline{\rm co}\{\sum\alpha_i\xi_i:\xi_i \in \partial h_i(x), \alpha \in \partial g(h(x))\}, $$ where $\overline{\rm co}$ denotes the weak-closed convex hull.
	\end{proposition}

\end{document}